\documentclass[11pt,letterpaper]{article}

\usepackage[T1]{fontenc}
\usepackage{lmodern}
\usepackage{microtype}
\usepackage{amsmath,amssymb,amsthm,mathtools}
\usepackage{xcolor}
\usepackage[margin=1.08in]{geometry}
\usepackage{enumitem}
\usepackage[
colorlinks=true,
linkcolor=magenta,
citecolor=cyan,
urlcolor=magenta,
pdftitle={Yang--Mills--Higgs Limits of Reflexive Higgs Sheaves},
pdfauthor={Tianshu Jiang and Yijian Zhang}
]{hyperref}
\hypersetup{linktoc=all,bookmarksnumbered=true,bookmarksopen=true}
\setlist[enumerate]{leftmargin=2.2em,itemsep=2pt,topsep=4pt}
\setlist[itemize]{leftmargin=2em,itemsep=2pt,topsep=4pt}
\allowdisplaybreaks
\numberwithin{equation}{section}
\theoremstyle{plain}
\newtheorem{theorem}{Theorem}[section]
\newtheorem{proposition}[theorem]{Proposition}
\newtheorem{lemma}[theorem]{Lemma}

\theoremstyle{definition}

\theoremstyle{remark}

\newcommand{\cE}{\mathcal E}
\newcommand{\cB}{\mathcal B}

\newcommand{\cI}{\mathcal I}
\newcommand{\Gr}{\operatorname{Gr}}
\newcommand{\Sing}{\operatorname{Sing}}
\newcommand{\Tor}{\operatorname{Tor}}
\newcommand{\Id}{\operatorname{Id}}
\newcommand{\Vol}{\operatorname{Vol}}
\newcommand{\rank}{\operatorname{rank}}
\newcommand{\im}{\operatorname{im}}
\newcommand{\Sat}{\operatorname{Sat}}
\newcommand{\tr}{\operatorname{tr}}
\newcommand{\ii}{\sqrt{-1}}
\newcommand{\ddV}{\,dV_{\omega}}
\newcommand{\HNS}{\mathrm{HNS}}

\title{Yang--Mills--Higgs Limits of Reflexive Higgs Sheaves}
\author{
Tianshu Jiang\setcounter{footnote}{1}%
\thanks{Nanjing University, Nanjing, China.
Email: \href{mailto:jts2026@nju.edu.cn}
{\nolinkurl{jts2026@nju.edu.cn}}.}
\and
Yijian Zhang\setcounter{footnote}{0}%
\thanks{Corresponding author.
University of Science and Technology of China, Hefei, China.
Email: \href{mailto:zyj_math@mail.ustc.edu.cn}
{\nolinkurl{zyj_math@mail.ustc.edu.cn}}.}
}
\date{}

\begin{document}
\maketitle

\begin{abstract}
For stable Higgs bundles, Hermitian--Einstein metrics can be obtained as
limits of the Yang--Mills--Higgs (Y--M--H) flow. We study the algebraic
structure of subsequential limits of this flow for reflexive Higgs sheaves
over compact Kähler manifolds, without assuming stability. We prove that
the limiting reflexive Higgs sheaf is isomorphic to the double dual of the
graded sheaf associated with the Harder--Narasimhan--Seshadri filtration of
the initial sheaf.
\end{abstract}

\noindent{\small
\par}

\pdfbookmark[1]{\contentsname}{contents}
\tableofcontents
\clearpage

\section{Introduction}

Let $(X,\omega)$ be a connected compact Kähler manifold, and let $\mathcal E$ be a coherent reflexive sheaf on $X$ endowed
with a Higgs field
$\theta\colon\mathcal E\to\mathcal E\otimes_{\mathcal O_X}\Omega_X^1$
satisfying $\theta\wedge\theta=0$.
We first recall the Harder--Narasimhan--Seshadri filtration of
$(\mathcal E,\theta)$. The $\omega$-degree and $\omega$-slope of a
torsion-free Higgs sheaf $(\mathcal F,\theta)$ are
\begin{equation}\label{degree and slope}
      \mathrm{deg}_\omega(\mathcal F) = \int_X c_1(\mathcal F)\wedge \frac{\omega^{n-1}}{(n-1)!},\quad \mu_\omega(\mathcal F) = \frac{\mathrm{deg}_\omega(\mathcal F)}{\mathrm{rank}(\mathcal F)},
\end{equation}
where $c_1(\mathcal F)$ is the first Chern class of $\mathcal F$. The Higgs sheaf $(\mathcal F,\theta)$ is $\omega$-stable (resp. $\omega$-semistable) if for every proper $\theta$-invariant subsheaf $\mathcal F'\subset \mathcal F$ with $
0<\operatorname{rank}\mathcal F'<\operatorname{rank}\mathcal F
$, $\mu_\omega(\mathcal F') < \mu_\omega(\mathcal F)$ (resp. $\mu_\omega(\mathcal F')\leq \mu_\omega(\mathcal F)$).  The Harder--Narasimhan filtration of $(\mathcal F,\theta)$ consists of a strictly increasing sequence
\begin{equation*}
  0 = \mathcal F_0 \subsetneq \mathcal F_1 \subsetneq \cdots \subsetneq \mathcal F_k = \mathcal F
\end{equation*}
of $\theta$-invariant saturated subsheaves, such that the successive quotients $\mathcal F_s/\mathcal F_{s-1}$ are semistable and $\mu_\omega(\mathcal F_1/ \mathcal F_0) > \mu_\omega(\mathcal F_2/ \mathcal F_1) > \cdots > \mu_\omega(\mathcal F_k/ \mathcal F_{k-1})$. Secondly, a semistable Higgs sheaf $(\mathcal V,\theta)$ admits a Seshadri
(or Jordan--Hölder) filtration
\begin{equation*}
  0 = \mathcal V_0 \subsetneq \mathcal V_1 \subsetneq \cdots \subsetneq \mathcal V_l = \mathcal V
\end{equation*}
of $\theta$-invariant saturated subsheaves, such that the successive quotients $\mathcal V_s/\mathcal V_{s-1}$ are stable and
\[
  \mu_\omega(\mathcal V_s/\mathcal V_{s-1}) = \mu_\omega(\mathcal V).
\]
Refining the Harder--Narasimhan filtration of $(\mathcal E,\theta)$ by
Seshadri filtrations gives the Harder--Narasimhan--Seshadri filtration
\begin{equation}\label{HNS-filtration}
  0 = \mathcal E_0 \subsetneq \mathcal E_1 \subsetneq \cdots \subsetneq \mathcal E_m = \mathcal E
\end{equation}
with stable successive quotients $\mathcal E_s/\mathcal E_{s-1}$.
We denote the double dual of its graded sheaf, equipped with the induced
Higgs field, by $\Gr^{\HNS}(\mathcal E,\theta)^{\vee\vee}$.
This reflexive Higgs sheaf is unique up to Higgs isomorphism.

Let $\Sigma_{\mathcal E}$ be the non-locally-free locus of $\mathcal E$ and put
$X_0=X\setminus\Sigma_{\mathcal E}$. The set $\Sigma_{\mathcal E}$ is
analytic of complex codimension at least three
\cite[Corollary~5.5.20]{Kobayashi87}.
On $X_0$, the sheaf $\mathcal E$ is locally free. Let $h$ be a smooth
Hermitian metric on $\mathcal E|_{X_0}$.
Let $\Theta_h$ denote the curvature of its Chern connection, and let
$\Lambda_\omega$ be the adjoint of the Lefschetz operator
$\alpha\mapsto\omega\wedge\alpha$ with respect to the metric induced by
$\omega$.
We call $h$ \emph{admissible} \cite{BS94} if
\[
 \int_{X_0}|\Theta_h|_{h,\omega}^2\frac{\omega^n}{n!}<\infty,
 \quad
 \sup_{X_0}|\Lambda_\omega\Theta_h|_h<\infty.
\]
Let $\Theta^{\mathrm{HS}}_{h,\theta}$ denote the Hitchin--Simpson curvature
defined in Section~\ref{sec:heat-kernel-estimates}.

Consider the heat flow on $\mathcal E|_{X_0}$:
\begin{equation}\label{heat flow}
  h^{-1}\frac{\partial h}{\partial t}
  =-2\bigl(\ii\Lambda_\omega\Theta^{\mathrm{HS}}_{h,\theta}
             -\lambda\Id_{\mathcal E}\bigr),
  \quad h(0)=h_0,
\end{equation}
where $\lambda=2\pi\mu_\omega(\mathcal E)/\Vol(X,\omega)$.

Under a suitable global complex gauge transformation, the heat flow is equivalent to the following Yang--Mills--Higgs flow concerning the evolution of a Higgs pair \((A(t),\theta(t))\) where \(A(t)\) is a flow of \(h_0\)-unitary connections and \(\theta(t)\) is a flow of Higgs fields associated to it:
\begin{equation}\label{Yang-Mills-Higgs flow}
  \frac{\partial A}{\partial t} = -\nabla_A^* \Theta_A - (\partial_A- \bar \partial_A) \Lambda_\omega \ii[\theta,\theta^{*h_0}], \ \frac{\partial \theta}{\partial t} = -[\Lambda_\omega \ii(\Theta_A + [\theta,\theta^{*h_0}]),\theta].
\end{equation}

If $(\mathcal E,\theta)$ is stable, then the heat flow \eqref{heat flow} converges to the Hermitian--Einstein metric:
\[
\ii\Lambda_\omega\Theta^{\mathrm{HS}}_{h,\theta}=\lambda\operatorname{Id}_{\mathcal E}.
\]
This result was first established by Donaldson \cite{Donaldson85,Donaldson87} for ordinary holomorphic bundles over projective manifolds, and later by Simpson \cite{Simpson88} for Higgs bundles over general Kähler manifolds. For further developments we refer to \cite{BS94,LZ15Approx,LZZ16}. We remark that the work of Uhlenbeck and Yau \cite{UY86} is crucial, although they employed the continuity method to prove the existence of Hermitian--Einstein metrics.

If $(\mathcal E,\theta)$ is not stable, the heat flow \eqref{heat flow} or its equivalent Y--M--H flow \eqref{Yang-Mills-Higgs flow} does not converge in general.
However, for a suitable choice of initial metric and local unitary gauges, small-energy regularity and Uhlenbeck's local gauge theorem give a subsequential limit \((A_\infty, \theta_\infty)\) of the Y--M--H flow on an open dense subset of $X_0$ \cite{HongTian04}. It can be shown that the pair \((\bar{\partial}_{A_\infty}=A_\infty^{0,1}, \theta_\infty)\) gives us a limiting Higgs bundle. Bando--Siu's extension theorem
\cite[Theorem 2]{BS94} extends the limiting bundle to a reflexive sheaf
$\mathcal E_\infty$ on $X$; the limiting Higgs field also extends, giving us a Higgs sheaf
$(\mathcal E_\infty,\theta_\infty)$. Bando and Siu \cite{BS94} outlined
this construction for reflexive sheaves, and further details were given
in \cite{LZZ17Survey}. 

The following question was first raised informally by Atiyah and Bott \cite{AtiyahBott83} for bundles over Riemann surfaces: they expected the Harder--Narasimhan type to determine the limit of the Yang--Mills flow. Bando and Siu \cite{BS94} later conjectured that for a reflexive sheaf, the double dual of the graded sheaf associated with its HNS filtration should be isomorphic to the limiting sheaf determined by the flow. For vector bundles, the result was proved in complex dimension one \cite{Daskalopoulos92}, dimension two \cite{DW04}, and higher dimensions \cite{Sibley15}. The conjecture was finally resolved in \cite{LZZ18}.

For Higgs bundles, the corresponding result was proved in complex dimension one \cite{Wilkin08}, two \cite{LZ11} and in higher dimensions \cite{LZ17}. For reflexive Higgs sheaves, \cite{LZZ17Survey} treats the case in which the HNS filtration contains exactly one nontrivial subsheaf. More recent work addresses the analogous problem for twisted Higgs pairs on bundles \cite{PSZ24}. For generalizations to non-Kähler manifolds, we refer to \cite{NieZhang20,CLZZ26}.

Our main theorem proves this identification for reflexive Higgs sheaves with HNS filtrations of arbitrary length:

\begin{theorem}\label{thm:main-theorem}
$(\mathcal E_{\infty}, \theta_{\infty})$ is isomorphic to $\Gr^{\HNS}(\mathcal E,\theta)^{\vee\vee}$ as Higgs sheaves.
\end{theorem}

\subsection*{Organization of the paper}

We begin in Section~\ref{sec:regularization} by regularizing the reflexive
Higgs sheaf and its HNS filtration. Uniform heat kernel estimates and
bounds on the Higgs field then give the heat flow on $X_0$ and its
admissibility at positive times in Section~\ref{sec:heat-kernel-estimates}.
The energy identity and the decay of the dissipation are proved in
Section~\ref{sec:cutoff-dissipation}.

In Section~\ref{sec:slice-compactness}, we write the flow in gauge form
and use small energy regularity to obtain a subsequential limit in local
unitary gauges. We extend this limit to a reflexive Higgs sheaf and
decompose it into polystable spectral summands in
Section~\ref{sec:extension-polystability}.

Section~\ref{sec:hn-vector} identifies the Harder--Narasimhan type with
the spectral vector of the limit. Section~\ref{sec:residual-exhaustion}
then proves the isomorphism by successively splitting off the stable HNS
factors.

\section*{Acknowledgements}
\addcontentsline{toc}{section}{Acknowledgements}

We are very grateful to Professor Jiayu Li and Professor Xi Zhang for their careful guidance and valuable discussions. We also thank Professor Chuanjing Zhang for her helpful advice.

\section{Regularization of the Higgs sheaf and its HNS filtration}
\label{sec:regularization}

The resolution results for coherent sheaves and their filtrations recalled
below are consequences of Hironaka's work on resolution of singularities
\cite{Hironaka64} and flattening in complex-analytic geometry
\cite{Hironaka75}, although we don't use his profound theorems directly in this section.

For a coherent sheaf $\mathcal F$ on a complex manifold,
$\Tor(\mathcal F)$ denotes its torsion subsheaf. If
$\mathcal A\subset\mathcal C$ is a coherent subsheaf of a vector bundle,
its saturation is
\[
 \Sat_{\mathcal C}(\mathcal A)
 =\ker\bigl(\mathcal C\longrightarrow
       (\mathcal C/\mathcal A)/\Tor(\mathcal C/\mathcal A)\bigr).
\]

\begin{proposition}\label{prop:relative-sheaf-resolution}
Let $M$ be a connected compact Kähler manifold and $\mathcal F$ a
coherent sheaf of positive rank on $M$. Let $W\subset M$ be its locally
free locus. There is a modification $p:M'\to M$, obtained by finitely
many blowups along smooth centers lying over $M\setminus W$, with
$M'$ Kähler, such that $p$ is biholomorphic
over $W$ and
$p^*\mathcal F/\Tor(p^*\mathcal F)$ is locally free.
\end{proposition}

\begin{proof}
Grivaux's theorem
\cite[Theorem~4.11, p.~25 of the arXiv version]{Grivaux10}
states that a coherent sheaf $\mathcal F$ on a compact complex manifold
$M$ admits a finite composition of blowups along smooth centers
$p_0:M_0\to M$ such that $p_0^*\mathcal F$ has a locally free quotient
of rank $\rank\mathcal F$.
The statement does not explicitly assert that $p_0$ is biholomorphic
over the locally free locus $W$. We show that the modification can be
chosen to satisfy this additional property.

Put $r=\rank\mathcal F$ and $\mathcal F(x)=\mathcal F_x/\mathfrak m_x\mathcal F_x$.
Grivaux constructs the reduced analytic space of rank-$r$ fiber quotients
\[
 q:G(\mathcal F)=\coprod_{x\in M}
       \operatorname{Gr}_{\mathrm{quot}}(r,\mathcal F(x))\longrightarrow M.
\]
For a local presentation
$\mathcal O_H^a\xrightarrow{A}\mathcal O_H^b
\to\mathcal F|_H\to0$, its local model is
\[
 q^{-1}(H)=\bigl\{(x,K)\in H\times\Gr(b-r,b):
                         \im A(x)\subset K\bigr\}_{\mathrm{red}}.
\]
These models glue independently of the presentation. They show that $q$
is proper and that $q^{-1}(W)\simeq W$. By
\cite[Lemma~4.12, pp.~27--28 of the arXiv version]{Grivaux10}, there is
a tautological locally free quotient
$q^*\mathcal F\twoheadrightarrow\mathcal Q$ of rank $r$.

Let $G_0$ be the reduced closure of $q^{-1}(W)$ in $G(\mathcal F)$,
the component dominating $M$.
Then $q_0=q|_{G_0}$ is a proper modification. Take a resolution
$d:Z\to G_0$ that is an isomorphism over the regular locus
\cite[Section~1.2.1]{AKMW02}.
The map $h=q_0\circ d$ is biholomorphic over $W$.
Elimination of the indeterminacies of $h^{-1}$
\cite[Sections~0.3(2) and~1.2.4]{AKMW02} gives a modification
$p:M'\to M$, biholomorphic over $W$, and a holomorphic map $g:M'\to Z$
with $h\circ g=p$. The centers defining $p$ lie over $M\setminus W$.
The manifold $M'$ is Kähler by
\cite[Section~3.3.3, Proposition~3.24, pp.~80--81]{Voisin02}.
Pulling back the tautological quotient gives
\[
 p^*\mathcal F\twoheadrightarrow g^*d^*(\mathcal Q|_{G_0}).
\]
Its kernel has rank zero and hence is torsion; conversely, every torsion
section maps to zero in the locally free target. Thus its kernel is
$\Tor(p^*\mathcal F)$, proving the assertion.
\end{proof}

We use the following relative form of Sibley's resolution of saturated filtrations
\cite[Section~4.1, Proposition~4.3]{Sibley15}. The choice of blowup centers
is also discussed, for the HNS filtration, at the beginning of the proof
of Theorem~5.12 in the same reference (p.~34 of the arXiv version).

\begin{proposition}\label{prop:relative-filtration-resolution}
Let $\mathcal C$ be a holomorphic vector bundle on a connected compact
Kähler manifold $M$, and let
\[
 0=\mathcal S_0\subsetneq\mathcal S_1\subsetneq\cdots
   \subsetneq\mathcal S_m=\mathcal C
\]
be a finite filtration by coherent saturated subsheaves. Suppose that
$U\subset M$ is a dense open set with analytic complement and that all
$\mathcal C/\mathcal S_i$ are locally free on $U$.
There is a modification $p:M'\to M$, obtained by finitely many blowups
along smooth centers avoiding the successive inverse images of $U$,
with $M'$ Kähler, biholomorphic over $U$,
such that
\[
 \mathcal T_i=\Sat_{p^*\mathcal C}
       \im(p^*\mathcal S_i\longrightarrow p^*\mathcal C)
\]
form a filtration by holomorphic subbundles of $p^*\mathcal C$.
All ambient and successive quotients are locally free, and the filtration
over $p^{-1}(U)$ is the pullback of the original one.
\end{proposition}

For a holomorphic map $f:N\to M$, the pullback of a Higgs field includes
the cotangent map. Thus, throughout the paper, $f^*\alpha$ as a Higgs
field means
\[
 (\Id\otimes df)\circ f^*\alpha,
 \quad df:f^*\Omega_M^1\longrightarrow\Omega_N^1.
\]

\begin{proposition}\label{prop:ambient-flag}
Let $(\cE,\theta)$ be a reflexive Higgs sheaf of positive rank on
$(X,\omega)$, and let
\[
 0=\cE_0\subsetneq\cE_1\subsetneq\cdots\subsetneq\cE_m=\cE
\]
be a finite filtration by coherent saturated $\theta$-invariant
subsheaves. Let $X_0$ be the locally free locus of $\cE$, and let
$X_f\subset X_0$ be the common locally free locus of the terms
$\cE_i$, the quotients $\cE/\cE_i$, and the successive quotients
$\cE_i/\cE_{i-1}$.
There are Kähler manifolds $Y$ and $\widetilde Y$ and finite
compositions of blowups along smooth centers
\[
 \rho:Y\longrightarrow X,\quad
 \sigma:\widetilde Y\longrightarrow Y
\]
with the following properties.

The centers defining $\rho$ lie over $X\setminus X_0$, so $\rho$ is
biholomorphic over $X_0$. Set
\[
 \mathcal M=\rho^*\cE,\quad T=\Tor(\mathcal M),\quad
 \cB=\mathcal M/T.
\]
The sheaf $\cB$ is locally free.
The field $\gamma=\rho^*\theta$ preserves $T$ and descends uniquely to
an integrable Higgs field $\beta$ on $\cB$.
The canonical quotient $q:\mathcal M\twoheadrightarrow\cB$ is a Higgs morphism
and is an isomorphism over $Y_0=\rho^{-1}(X_0)$.

Put
\[
 \mathcal H_i=\im\bigl(\rho^*\cE_i\longrightarrow\mathcal M
                              \xrightarrow{q}\cB\bigr),\quad
 \mathcal S_i=\Sat_{\cB}(\mathcal H_i),\quad U=\rho^{-1}(X_f).
\]
The centers defining $\sigma$ avoid the successive inverse images of $U$,
so $\sigma$ is biholomorphic over $U$.
Writing $\cB'=\sigma^*\cB$, $\delta=\sigma^*\beta$, and
$\pi=\rho\circ\sigma$, the sheaves
\[
 \widetilde\cE_i=\Sat_{\cB'}
       \im(\sigma^*\mathcal S_i\longrightarrow\cB')
\]
form a strict $\delta$-invariant holomorphic subbundle filtration
\[
 0=\widetilde\cE_0\subsetneq\widetilde\cE_1\subsetneq\cdots
   \subsetneq\widetilde\cE_m=\cB'.
\]
Every $\widetilde\cE_i$, $\cB'/\widetilde\cE_i$, and
$\widetilde\cE_i/\widetilde\cE_{i-1}$ is locally free and carries the
induced integrable Higgs field. Inside $\cB'$ we have
\begin{equation}\label{eq:direct-saturated-transform}
 \widetilde\cE_i=\Sat_{\cB'}\im\bigl(
   \pi^*\cE_i\longrightarrow\pi^*\cE
                  \xrightarrow{\sigma^*q}\cB'\bigr),
\end{equation}
where $\pi^*\cE\simeq\sigma^*\rho^*\cE$.

The set $X_f$ is dense, and $X\setminus X_f$ is closed analytic of
complex codimension at least two. The map $\pi$ is biholomorphic over
$X_f$; there the resulting filtration and all its quotient Higgs bundles
agree with the pullbacks of the original ones.
\end{proposition}

\begin{proof}
Apply Proposition~\ref{prop:relative-sheaf-resolution} to $\cE$ to
obtain $\rho$ and $\cB=\rho^*\cE/\Tor(\rho^*\cE)$.
In local coordinates $x_1,\ldots,x_n$ on $X$, write
$\theta=\sum_a\theta_a\otimes dx_a$. The relation
$\theta\wedge\theta=0$ means that $[\theta_a,\theta_b]=0$, and
\[
 \gamma=\sum_a\rho^*\theta_a\otimes d(x_a\circ\rho).
\]
Pullback respects compositions, so $\gamma\wedge\gamma=0$.
For a torsion germ $t\in T_y$, choose a nonzero local holomorphic
function $f$ with $ft=0$. Then
$f\gamma(t)=\gamma(ft)=0$.
Since $\Omega_Y^1$ is locally free, each coefficient of $\gamma(t)$ is
torsion, and $\gamma(T)\subset T\otimes\Omega_Y^1$.
Thus $\gamma$ descends uniquely to $\beta$, with
\[
 \beta q=(q\otimes\Id)\gamma,\quad
 (\beta\wedge\beta)q=(q\otimes\Id)(\gamma\wedge\gamma)=0.
\]
Surjectivity of $q$ proves that $\beta$ is integrable and that $q$ is a
Higgs morphism. Over $Y_0$, the sheaf $\mathcal M$ is locally free, so $T=0$ and
$q$ is an isomorphism.

We next check that saturation preserves Higgs invariance. Suppose
$\mathcal A\subset\mathcal C$ is invariant under a Higgs field $\alpha$
on a vector bundle, and put $\mathcal S=\Sat_{\mathcal C}(\mathcal A)$.
At a stalk, $s\in\mathcal S$ means that $fs\in\mathcal A$ for some
nonzero local function $f$.
The image of $\alpha(s)$ in
$(\mathcal C/\mathcal S)\otimes\Omega^1$ is therefore killed by $f$.
This target is torsion-free, so the image vanishes.
Images of Higgs morphisms are invariant as well, since tensoring with
the locally free sheaf $\Omega^1$ preserves images.
Consequently, $\mathcal H_i$ and $\mathcal S_i$ are $\beta$-invariant,
and the $\mathcal S_i$ form a saturated coherent filtration of $\cB$.
Its endpoints are $0$ and $\cB$, since $q$ is surjective.

The sheaves $\cE_i$ are torsion-free, and
$\cE_i/\cE_{i-1}$ embeds into the torsion-free sheaf
$\cE/\cE_{i-1}$. Every torsion-free coherent sheaf on a complex
manifold is locally free off an analytic subset of complex codimension
at least two. This proves the assertions about $X_f$.
On $X_f$, the sequences
\[
 0\longrightarrow\cE_i\longrightarrow\cE
        \longrightarrow\cE/\cE_i\longrightarrow0
\]
split locally. Over $U$, the map $q$ therefore identifies
$\mathcal H_i$ with the pullback subbundle of $\cE_i$.
Its ambient quotient is locally free, so saturation does not change it.
Thus all $\cB/\mathcal S_i$ are locally free on $U$.
Each nonzero successive quotient of the original strict filtration is
torsion-free and has positive rank. These ranks are preserved on $U$,
so the filtration by $\mathcal S_i$ is strict.

Apply Proposition~\ref{prop:relative-filtration-resolution} to this
filtration and $U$. It gives $\sigma$ and the subbundles
$\widetilde\cE_i$, with all ambient and successive quotients locally
free. The maps $\sigma^*\mathcal S_i\to\cB'$ are Higgs morphisms.
Their images and saturations are therefore invariant under
$\delta=\sigma^*\beta$.
The coordinate calculation above also gives $\delta\wedge\delta=0$.
Restricting this identity to invariant subsheaves and passing to
quotients proves integrability of all induced fields, because tensoring
with $\Omega^2_{\widetilde Y}$ is exact.

To prove \eqref{eq:direct-saturated-transform}, set
\[
 \begin{aligned}
 A_i&=\im(\sigma^*\mathcal H_i\longrightarrow\cB'),&
 J_i&=\im(\sigma^*\mathcal S_i\longrightarrow\cB'),\\
 D_i&=\mathcal S_i/\mathcal H_i
       =\Tor(\cB/\mathcal H_i).
 \end{aligned}
\]
Right exactness of pullback gives
\[
 \pi^*\cE_i\twoheadrightarrow\sigma^*\mathcal H_i,
 \quad
 \sigma^*\mathcal H_i\longrightarrow\sigma^*\mathcal S_i
                     \longrightarrow\sigma^*D_i\longrightarrow0.
\]
The first map identifies $A_i$ with the image in
\eqref{eq:direct-saturated-transform}; the second sequence gives a
surjection $\sigma^*D_i\twoheadrightarrow J_i/A_i$.
A local nonzero annihilator of $D_i$ stays nonzero after pullback by
the modification $\sigma$. Hence $\sigma^*D_i$, and therefore $J_i/A_i$,
is torsion. It follows that
\[
 \Sat_{\cB'}(A_i)=\Sat_{\cB'}(J_i)=\widetilde\cE_i,
\]
which proves \eqref{eq:direct-saturated-transform}.
Both sides carry the restriction of $\delta$, so this equality also
identifies their Higgs fields.

Finally, $\rho$ preserves $X_0$ and $\sigma$ preserves $U$, so $\pi$
is biholomorphic over $X_f$.
On this open set all defining bundle sequences split locally, and
neither saturation changes the pulled-back filtration.
The identity $d(\rho\circ\sigma)=d\sigma\circ\sigma^*d\rho$ identifies
the induced Higgs fields with the pullbacks of the original ones.
\end{proof}

For the remainder of the paper, we apply
Proposition~\ref{prop:ambient-flag} to the HNS filtration
\eqref{HNS-filtration} and fix the modifications constructed in its proof.

\section{Heat flow and admissibility}\label{sec:heat-kernel-estimates}

We recall the heat-flow construction and estimates of Li, Zhang, and Zhang
\cite[Sections~2 and~3]{LZZ16}, using the notation and regularization of
Section~\ref{sec:regularization}.

For a Hermitian metric $h$ on $\cE|_{X_0}$, let
$\partial_h+\bar\partial$ be the Chern connection, $\Theta_h$ its curvature,
and $\theta^{*h}$ the metric adjoint of the Higgs field.
The Hitchin--Simpson connection and its curvature are
\begin{align*}
 \nabla^{\mathrm{HS}}&=\partial_h+\bar\partial+\theta+\theta^{*h},\\
 \Theta^{\mathrm{HS}}_{h,\theta}&=(\nabla^{\mathrm{HS}})^2
 =\Theta_h+[\theta,\theta^{*h}]
   +\partial_h\theta+\bar\partial\theta^{*h},
\end{align*}
where the bracket denotes the graded commutator.  We write
\[
 \Theta^{1,1}_{h,\theta}=\Theta_h+[\theta,\theta^{*h}]
\]
for its $(1,1)$ component, so that
$\ii\Lambda_\omega\Theta^{\mathrm{HS}}_{h,\theta}
 =\ii\Lambda_\omega\Theta^{1,1}_{h,\theta}$.
Throughout, for a Kähler form $\kappa$ on a compact complex $n$-fold $M$, set
\[
 dV_{\kappa}:=\frac{\kappa^n}{n!},\quad
 \Vol(M,\kappa):=\int_M dV_{\kappa}.
\]
On scalar functions, we write
$\Delta_\kappa=2\ii\Lambda_\kappa\partial\bar\partial$.

Choose a Kähler form $\eta$ on $Y$ with
$\int_Y\eta^n/n!=1$, and a smooth Hermitian metric $\widehat h$ on $\cB$.
Using the bundle isomorphism over $Y_0=\rho^{-1}(X_0)$,
we define the initial metric $h_0$ appearing in the Introduction by
\[
 (\rho|_{Y_0})^*h_0=\widehat h|_{Y_0}.
\]
For $0<\epsilon\leq1$, put
\[
 \omega_{\epsilon}=\rho^*\omega+\epsilon\eta.
\]
Define
\[
 \lambda_{\epsilon}
 =\frac{2\pi\mu_{\omega_{\epsilon}}(\cB)}
        {\Vol(Y,\omega_{\epsilon})},
 \quad
 \lambda=\frac{2\pi\mu_{\omega}(\cE)}{\Vol(X,\omega)}.
\]
Let $h_{\epsilon}(t)$ solve
\begin{equation}\label{eq:regularized-flow}
 h_{\epsilon}^{-1}\partial_th_{\epsilon}
 =-2\Phi_{\epsilon},\quad
 \Phi_{\epsilon}
 =\ii\Lambda_{\omega_{\epsilon}}
 \bigl(\Theta_{h_{\epsilon}}+[\beta,\beta^{*h_{\epsilon}}]\bigr)
 -\lambda_{\epsilon}\Id,
\end{equation}
with $h_{\epsilon}(0)=\widehat h$.  For each fixed $\epsilon>0$,
$(Y,\omega_\epsilon)$ is a compact Kähler manifold.  The long-time
existence result for the Higgs bundle heat flow \cite{Simpson88,LZZ16}
therefore gives a smooth solution of \eqref{eq:regularized-flow}
for all $t\geq0$.

Let $D=\{y\in Y:\rank_{\mathbb C}d\rho_y<n\}$ be the exceptional locus
of $\rho$, and put $V_\epsilon=\Vol(Y,\omega_\epsilon)$.
Since $\rho$ has degree one and
$\rho^*\omega\leq\omega_\epsilon\leq\omega_1$, we have
\[
 0<V_-:=\Vol(X,\omega)\leq V_\epsilon
 \leq V_+:=\Vol(Y,\omega_1).
\]
By \cite[(2.3)]{LZZ16}, we have
\begin{equation}\label{eq:initial-phi-bound}
 \sup_{0<\epsilon\leq1}\int_Y
 |\Phi_\epsilon(0)|_{\widehat h}\,dV_{\omega_\epsilon}<\infty.
\end{equation}

\begin{lemma}\label{lem:fixed-eta-analytic-gate}
Let $H_\epsilon(x,y,t)$ be the scalar heat kernel acting against
$dV_{\omega_\epsilon}$, normalized by
$\int_YH_\epsilon(x,y,t)dV_{\omega_\epsilon}(y)=1$.
There is a constant $C>0$, depending only on the fixed geometry
$(X,\omega)$, $\rho:Y\to X$, and $(Y,\eta)$, such that
\begin{equation}\label{eq:uniform-heat-kernel}
 H_\epsilon(x,y,t)\leq C\left(
 t^{-n}\exp\left(-\frac{d_{\omega_\epsilon}(x,y)^2}{10t}\right)+1
 \right),\quad t>0,\quad 0<\epsilon\leq1.
\end{equation}
\end{lemma}

\begin{proof}
See \cite[Proposition~2.1]{LZZ16}, with its parameter $\tau=6$.
\end{proof}

We also use the bound \cite[(2.11)]{LZZ16}
\begin{equation}\label{eq:regularized-phi-pointwise}
 \sup_Y|\Phi_\epsilon(t)|_{h_\epsilon(t)}\leq C(1+t^{-n}),
 \quad t>0,\quad 0<\epsilon\leq1,
\end{equation}
where $C$ is independent of $\epsilon$ and $t$.

\begin{lemma}\label{lem:initial-higgs-integrability}
There is a constant $\delta_0>0$ such that
\begin{equation}\label{eq:initial-analytic-inputs}
 \sup_{0<\epsilon\leq1}\int_Y
 |\beta|_{\widehat h,\omega_\epsilon}^{2(1+\delta_0)}
 dV_{\omega_\epsilon}<\infty.
\end{equation}
\end{lemma}

\begin{proof}
See \cite[(3.3)--(3.4)]{LZZ16}, using the integrability exponent supplied
by \cite[Lemma~4.4 and (4.10)--(4.12)]{LZZ18} for the successive blowups
defining $\rho$.
\end{proof}

\begin{proposition}\label{prop:positive-flow}
There is a subsequence $\epsilon_j\downarrow0$ such that
$h_{\epsilon_j}(t)$ converges to a metric $h(t)$ in $C^{\infty}_{\mathrm{loc}}(X_0\times[0,\infty))$, which satisfies
\begin{equation}\label{eq:singular-flow}
 h^{-1}\partial_th=-2\Phi,\quad h(0)=h_0,\quad
 \Phi=\ii\Lambda_{\omega}
 \Theta^{1,1}_{h,\theta}-\lambda\Id.
\end{equation}
For every $\tau>0$, there is a finite constant $C_\tau$, independent of
$t\geq\tau$, such that
\begin{equation}\label{eq:positive-time-bounds}
 \begin{aligned}
 \|\Phi(t)\|_{L^{\infty}}+\|\theta\|_{L^{\infty}(h(t),\omega)}^2
   +\|\Theta^{\mathrm{HS}}_{h(t),\theta}\|_{L^2}^2
  +\|\Theta_{h(t)}\|_{L^2}^2
   +\|\Lambda_{\omega}\Theta_{h(t)}\|_{L^{\infty}}\leq C_{\tau}.
 \end{aligned}
\end{equation}
Thus $h(t)$ is admissible for every $t>0$.
\end{proposition}

\begin{proof}
See \cite[Lemma~2.6]{LZZ16} for the convergence and the limiting flow;
we fix the resulting subsequence once and for all.
For the pointwise bounds, see \cite[(2.11) and Proposition~3.4]{LZZ16}.
In the proof of the Higgs bound, we take
$0<2a<\min\{\delta_0,1/4\}$, using
Lemma~\ref{lem:initial-higgs-integrability}; the estimates in
\cite[(3.15)--(3.22)]{LZZ16} hold for this choice as well.
For $n\geq2$, see \cite[(3.30) and Corollary~3.5]{LZZ16} for the
curvature bounds and admissibility, uniformly for $t\geq\tau$.
When $n=1$, contraction with $\omega$ controls $(1,1)$-forms pointwise
and $\Theta^{\mathrm{HS}}_{h,\theta}=\Theta^{1,1}_{h,\theta}$, giving
the same curvature bounds.
\end{proof}

\section{Energy identity and decay estimates}\label{sec:cutoff-dissipation}

Let $\nabla_{h(t)}$ denote the Chern connection on $\cE|_{X_0}$.
For $s\in\Gamma^\infty(\operatorname{End}(\cE|_{X_0}))$, the induced
Hitchin--Simpson connection acts by
\[
 \nabla^{\mathrm{HS}}(t)s
 =\nabla_{h(t)}s+[\theta+\theta^{*h(t)},s].
\]
Define the dissipation along $h(t)$ by
\begin{equation}\label{eq:dissipation-definition}
 \cI(t)=\int_{X_0}
 \bigl|\nabla^{\mathrm{HS}}(t)\Phi(t)\bigr|^2\ddV
 =\int_{X_0}\left(
 |\nabla_{h(t)}\Phi(t)|^2+2|[\Phi(t),\theta]|^2\right)\ddV.
\end{equation}
All norms here are taken with respect to $h(t)$ and $\omega$.
Set
\[
 E(t)=\int_{X_0}|\Phi(t)|^2\ddV,\quad t>0.
\]

Let $\nabla_{h_\epsilon(t)}$ be the Chern connection on $\cB$.
For $s\in\Gamma^\infty(\operatorname{End}\cB)$, the induced
Hitchin--Simpson connection acts by
\[
 \nabla^{\mathrm{HS}}_\epsilon(t)s
 =\nabla_{h_\epsilon(t)}s+
   [\beta+\beta^{*h_\epsilon(t)},s].
\]
For the regularized flow put
\begin{align*}
 \cI_{\epsilon}(t)
 &=\int_Y
   \bigl|\nabla^{\mathrm{HS}}_\epsilon(t)
                 \Phi_{\epsilon}(t)\bigr|^2
       dV_{\omega_{\epsilon}},\\
 E_{\epsilon}(t)
 &=\int_Y|\Phi_{\epsilon}(t)|^2
       dV_{\omega_{\epsilon}}.
\end{align*}

By the estimates in Section~\ref{sec:heat-kernel-estimates}, for every
$\tau>0$ we can choose $M_\tau<\infty$ such that
\[
 \begin{aligned}
 |\Phi_\epsilon|&\leq M_\tau
 &&\text{on }Y\times[\tau,\infty),\quad 0<\epsilon\leq1,\\
 |\Phi|&\leq M_\tau
 &&\text{on }X_0\times[\tau,\infty).
 \end{aligned}
\]

\begin{lemma}\label{lem:energy-identities}
For every $0<\epsilon\leq1$ and $0\leq a<b<\infty$,
\begin{equation}\label{eq:endpoint-identity-epsilon}
 E_{\epsilon}(b)-E_{\epsilon}(a)
 +2\int_a^b\cI_{\epsilon}(t)\,dt=0.
\end{equation}
For every $0<t_1<t_2<\infty$,
\begin{equation}\label{eq:endpoint-identity-limit}
 E(t_2)-E(t_1)+2\int_{t_1}^{t_2}\cI(t)\,dt=0.
\end{equation}
Consequently, for every $\tau>0$,
\begin{equation}\label{eq:total-dissipation}
 \int_\tau^\infty\cI(t)\,dt\leq\frac12E(\tau)<\infty.
\end{equation}
\end{lemma}

\begin{proof}
Fix $\tau>0$.  The endomorphism $\Phi_\epsilon$ is self-adjoint.  On a real orthonormal tangent
frame, $\nabla_{h_\epsilon}\Phi_\epsilon$ takes values in
self-adjoint endomorphisms, whereas
$[\beta+\beta^{*h_\epsilon},\Phi_\epsilon]$ takes values in
skew-adjoint endomorphisms.  Thus the mixed term in the expansion of
$|\nabla^{\mathrm{HS}}_\epsilon(t)\Phi_\epsilon|^2$ vanishes:
\[
 2\operatorname{Re}\left\langle
 \nabla_{h_\epsilon}\Phi_\epsilon,\,
 [\beta+\beta^{*h_\epsilon},\Phi_\epsilon]
 \right\rangle_{h_\epsilon,\omega_\epsilon}=0.
\]
Since $\Phi_\epsilon$ is self-adjoint and the Chern connection preserves
adjoints, we have
\[
 (\partial_{h_\epsilon}\Phi_\epsilon)^{*h_\epsilon}
 =\bar\partial_{\cB}\Phi_\epsilon,\quad
 [\beta,\Phi_\epsilon]^{*h_\epsilon}
 =-[\beta^{*h_\epsilon},\Phi_\epsilon].
\]
Taking adjoints preserves norms, and $(1,0)$-forms are orthogonal to
$(0,1)$-forms.  Hence
\[
 \bigl|\nabla^{\mathrm{HS}}_\epsilon(t)\Phi_\epsilon\bigr|^2
 =|\nabla_{h_\epsilon}\Phi_\epsilon|^2
   +2|[\Phi_\epsilon,\beta]|^2
 =2\bigl|\bar\partial_{\cB}\Phi_\epsilon
          +[\beta,\Phi_\epsilon]\bigr|^2.
\]
Along \eqref{eq:regularized-flow}, we have
\cite[(2.4)--(2.6)]{LZZ16}
\begin{equation}\label{eq:regularized-scalar-identity}
 (\Delta_{\omega_{\epsilon}}-\partial_t)
 |\Phi_{\epsilon}|^2
 =2\bigl|\nabla^{\mathrm{HS}}_\epsilon(t)
              \Phi_\epsilon\bigr|^2.
\end{equation}

Fix $0\leq a<b<\infty$.  Integration over the compact regularization gives
\eqref{eq:endpoint-identity-epsilon}.

By Proposition~\ref{prop:positive-flow}, $h_{\epsilon_j}\to h$ in
$C^\infty_{\mathrm{loc}}(X_0\times(0,\infty))$.
Letting $j\to\infty$ in \eqref{eq:regularized-scalar-identity} gives
\begin{equation}\label{eq:scalar-identity-limit}
 (\Delta_{\omega}-\partial_t)|\Phi|^2
 =2\bigl|\nabla^{\mathrm{HS}}(t)\Phi\bigr|^2.
\end{equation}
The same self-adjointness argument for the metric $h(t)$ gives
\[
 \bigl|\nabla^{\mathrm{HS}}(t)\Phi\bigr|^2
 =|\nabla_{h(t)}\Phi|^2+2|[\Phi,\theta]|^2
 =2|\bar\partial_{\cE}\Phi+[\theta,\Phi]|^2.
\]
Reflexivity gives $\operatorname{codim}_{\mathbb C}\Sigma_{\mathcal E}\geq3$.
By \cite[(3.14), (3.17)]{LZZ16}, there exist
$\chi_\nu\in C_c^\infty(X_0)$ with
$0\leq\chi_\nu\uparrow1$ pointwise on $X_0$ and
\begin{equation}\label{eq:capacity-cutoff-norms}
 \int_X|d\chi_\nu|^4\ddV\longrightarrow0,
 \quad
 \kappa_\nu:=\int_X|d\chi_\nu|^2\ddV\longrightarrow0.
\end{equation}
For $t\geq\tau$, let
\[
 E_\nu(t)=\int_{X_0}\chi_\nu^2|\Phi|^2\ddV,
 \quad
 \cI_\nu(t)=\int_{X_0}\chi_\nu^2
       \bigl|\nabla^{\mathrm{HS}}(t)\Phi\bigr|^2\ddV.
\]
Using $d|\Phi|^2=2\operatorname{Re}\langle\nabla_h\Phi,\Phi\rangle$,
we multiply \eqref{eq:scalar-identity-limit} by $\chi_\nu^2$ and integrate
by parts to obtain
\begin{align}
 B_\nu(t)&=-4\operatorname{Re}\int_{X_0}\chi_\nu
 \langle\nabla_h\Phi,(d\chi_\nu)\otimes\Phi\rangle\ddV.
 \label{eq:cutoff-boundary}\\
 \frac{d}{dt}E_\nu(t)&=B_\nu(t)-2\cI_\nu(t).
 \label{eq:cutoff-differential}
\end{align}
The Cauchy--Schwarz inequality and Young's inequality give
\begin{equation}\label{eq:cutoff-young}
 |B_\nu(t)|\leq\cI_\nu(t)+4M_\tau^2\kappa_\nu.
\end{equation}
Define
\[
 C_\nu(\tau)=M_\tau^2\Vol(X,\omega)+8M_\tau^2\kappa_\nu.
\]
Consequently, whenever $\tau\leq t_1<t_2\leq t_1+2$,
\begin{equation}\label{eq:cutoff-window-bound}
 \int_{t_1}^{t_2}\cI_\nu(t)\,dt
 \leq C_\nu(\tau).
\end{equation}
By \eqref{eq:cutoff-window-bound} and the Cauchy--Schwarz inequality,
\begin{equation}\label{eq:cutoff-remainder}
 \int_{t_1}^{t_2}|B_\nu(t)|\,dt
 \leq4C_\nu(\tau)^{1/2}(2M_\tau^2\kappa_\nu)^{1/2}
 \longrightarrow0.
\end{equation}

Monotone convergence at the two endpoints gives
\[
 E_\nu(t_a)\uparrow E(t_a),\quad a=1,2.
\]
For the nonnegative cutoff dissipation density, monotone convergence gives
\[
 \int_{t_1}^{t_2}\cI_\nu(t)\,dt
 \uparrow\int_{t_1}^{t_2}\cI(t)\,dt<\infty,
\]
where finiteness follows from \eqref{eq:cutoff-window-bound}.
By \eqref{eq:cutoff-remainder},
$\int_{t_1}^{t_2}B_\nu(t)\,dt\to0$.
Letting
$\nu\to\infty$ in the integrated identity yields
\eqref{eq:endpoint-identity-limit}.
Partitioning extends the identity from windows of length at most two to
every finite window in $(0,\infty)$, since $\tau>0$ is arbitrary.
Since $E\geq0$, it also gives \eqref{eq:total-dissipation}.
\end{proof}

\begin{lemma}
\label{lem:regularized-dissipation-growth}
For every $\tau>0$, $0<\epsilon\leq1$, and
$\tau\leq s\leq t\leq s+1$,
\begin{equation}\label{eq:one-sided-growth}
 \cI_{\epsilon}(t)
 \leq e^{4M_{\tau}(t-s)}\cI_{\epsilon}(s).
\end{equation}
Here $M_\tau$ is the constant fixed above, satisfying
$|\Phi_\epsilon|\leq M_\tau$ on $Y\times[\tau,\infty)$
for all $0<\epsilon\leq1$.
\end{lemma}

\begin{proof}
Fix $\tau>0$ and $0<\epsilon\leq1$.  Write the Hitchin--Simpson
connection on $\cB$ as $D''+D'_{h_\epsilon}$, where
\[
 D''=\bar\partial+\beta,\quad
 D'_{h_\epsilon}=\partial_{h_\epsilon}+\beta^{*h_\epsilon}.
\]
We use the same notation for the induced operators on
$\operatorname{End}(\cB)$-valued forms.  Set
\[
 \Box''_{\epsilon,t}
 :=D''(D'')^{*_{h_\epsilon(t),\omega_\epsilon}}
   +(D'')^{*_{h_\epsilon(t),\omega_\epsilon}}D''.
\]
Adjoints are taken with respect to $h_\epsilon(t)$ and $\omega_\epsilon$.
The operator $D''$ is independent of time, while its adjoint depends on $t$.
On zero-forms, $\Box''_{\epsilon,t}=(D'')^{*_{h_\epsilon(t),\omega_\epsilon}}D''$.

We first derive the evolution equation for $\Phi_\epsilon$.  Put
\[
 s_\epsilon
 =h_\epsilon^{-1}\partial_th_\epsilon=-2\Phi_\epsilon.
\]
The Kähler form $\omega_\epsilon$ is fixed in time.  Variation of the
Chern curvature and of the metric adjoint of the Higgs field gives
\[
 \partial_t\Theta_{h_\epsilon}
 =\bar\partial_{\cB}(\partial_{h_\epsilon}s_\epsilon),
 \quad
 \partial_t(\beta^{*h_\epsilon})
 =[\beta^{*h_\epsilon},s_\epsilon].
\]
Since $\lambda_\epsilon$ is independent of $t$, these identities imply
\begin{align}
 \partial_t\Phi_\epsilon
 &=\ii\Lambda_{\omega_\epsilon}\left(
   \bar\partial_{\cB}(\partial_{h_\epsilon}s_\epsilon)
   +[\beta,[\beta^{*h_\epsilon},s_\epsilon]]\right)\nonumber\\
 &=\ii\Lambda_{\omega_\epsilon}
       D''D'_{h_\epsilon}s_\epsilon.
 \label{eq:regularized-phi-variation}
\end{align}
In the second line, the $(2,0)$ and $(0,2)$ components of
$D''D'_{h_\epsilon}s_\epsilon$ disappear under
$\Lambda_{\omega_\epsilon}$.

The Higgs Kähler identity
\[
 (D'')^{*_{h_\epsilon,\omega_\epsilon}}
 =-\ii[\Lambda_{\omega_\epsilon},D'_{h_\epsilon}]
\]
shows that, for $u\in\Gamma^\infty(\operatorname{End}\cB)$,
\begin{align*}
 (D'')^{*_{h_\epsilon,\omega_\epsilon}}D''u
 &=
 \ii\Lambda_{\omega_\epsilon}D''D'_{h_\epsilon}u
 -\ii\Lambda_{\omega_\epsilon}
       [D'_{h_\epsilon},D'']u,\\
 \ii\Lambda_{\omega_\epsilon}
       [D'_{h_\epsilon},D'']u
 &=[\Phi_\epsilon+\lambda_\epsilon\Id,u].
\end{align*}
For details of the Kähler identities and curvature calculation, see
\cite[Lemma~3.1(a) and the proof of part~(c), p.~876]{Simpson88}.
For $u=\Phi_\epsilon$, the last commutator vanishes.  Combining this
observation with \eqref{eq:regularized-phi-variation} and
$s_\epsilon=-2\Phi_\epsilon$ gives
\begin{equation}\label{eq:regularized-phi-evolution}
 \partial_t\Phi_\epsilon
 =-2\Box''_{\epsilon,t}\Phi_\epsilon.
\end{equation}

The endomorphism norm also depends on time.  Set
\[
 J_\epsilon(t)=\int_Y|D''\Phi_\epsilon|^2
       dV_{\omega_\epsilon}.
\]
At a fixed point and time, choose an
$h_\epsilon$-orthonormal frame and an
$\omega_\epsilon$-unitary coframe, and write an
endomorphism-valued form as $v=\sum_\alpha v_\alpha\otimes e^\alpha$.
Because
\[
 \partial_t(v_\alpha^{*h_\epsilon})
 =(\partial_tv_\alpha)^{*h_\epsilon}
   +[v_\alpha^{*h_\epsilon},s_\epsilon],
\]
one obtains
\[
 \partial_t|v|^2
 -2\operatorname{Re}\langle\partial_tv,v\rangle
 =
 \sum_\alpha\left\{
 \tr(s_\epsilon v_\alpha v_\alpha^{*h_\epsilon})
 -\tr(s_\epsilon v_\alpha^{*h_\epsilon}v_\alpha)\right\}.
\]
The two trace terms satisfy
\[
 \left|
 \partial_t|v|^2
 -2\operatorname{Re}\langle\partial_tv,v\rangle
 \right|
 \leq2\|s_\epsilon\|_{\mathrm{op}}|v|^2.
\]
Apply this estimate to $v=D''\Phi_\epsilon$.  Since
$\|s_\epsilon\|_{\mathrm{op}}
 \leq2|\Phi_\epsilon|\leq2M_\tau$ for $t\geq\tau$,
\begin{align*}
 J_\epsilon'(t)
 &\leq
 2\operatorname{Re}\int_Y
 \langle D''(\partial_t\Phi_\epsilon),
         D''\Phi_\epsilon\rangle\,dV_{\omega_\epsilon}
 +4M_\tau J_\epsilon(t)\\
 &=-4\operatorname{Re}\int_Y
 \langle D''\Box''_{\epsilon,t}\Phi_\epsilon,
         D''\Phi_\epsilon\rangle\,dV_{\omega_\epsilon}
 +4M_\tau J_\epsilon(t)\\
 &=-4\int_Y
 |\Box''_{\epsilon,t}\Phi_\epsilon|^2
 dV_{\omega_\epsilon}
 +4M_\tau J_\epsilon(t)\\
 &\leq4M_\tau J_\epsilon(t).
\end{align*}

Finally, self-adjointness of $\Phi_\epsilon$ and orthogonality of types
give
\[
 2|D''\Phi_\epsilon|^2
 =|\nabla_{h_\epsilon}\Phi_\epsilon|^2
   +2|[\Phi_\epsilon,\beta]|^2
 =\bigl|\nabla^{\mathrm{HS}}_\epsilon(t)
             \Phi_\epsilon\bigr|^2.
\]
Since $\cI_\epsilon=2J_\epsilon$, the preceding inequality gives
\[
 \cI_\epsilon'(t)\leq4M_\tau\cI_\epsilon(t),\quad t\geq\tau.
\]
Integrating yields
\eqref{eq:one-sided-growth}.
The coefficient depends only on the bound $M_\tau$ for $\Phi_\epsilon$,
not on a bound for $\beta$ or the curvature of $\omega_\epsilon$.
\end{proof}

\begin{proposition}\label{prop:dissipation}
For every $\tau>0$, the function $\cI$ is integrable on $[\tau,\infty)$ and, for every
$t\geq\tau+1$,
\begin{equation}\label{eq:dissipation-decay}
 \cI(t)\leq e^{4M_{\tau}}
 \int_{t-1}^{t}\cI(s)\,ds\longrightarrow0
 \quad (t\to\infty).
\end{equation}
\end{proposition}

\begin{proof}
Integrability follows from Lemma~\ref{lem:energy-identities}.
For each fixed pair $0<\tau\leq t$, the uniform bound
$|\Phi_\epsilon(t)|\leq M_\tau$,
local smooth convergence off the exceptional locus $D$, and
$\omega_\epsilon\leq\omega_1$ show, by removing a neighborhood of $D$
of arbitrarily small $\omega_1$-volume, that
\begin{equation}\label{eq:endpoint-convergence}
 E_{\epsilon_j}(t)\longrightarrow E(t).
\end{equation}
For every finite window $[a,b]\subset(0,\infty)$, choose
$0<\tau\leq a$.  The two exact endpoint identities
\eqref{eq:endpoint-identity-epsilon} and
\eqref{eq:endpoint-identity-limit} therefore imply
\begin{equation}\label{eq:window-convergence}
 \int_a^b\cI_{\epsilon_j}(t)\,dt
 \longrightarrow\int_a^b\cI(t)\,dt.
\end{equation}
Finally, choose compact sets $K_N\Subset X_0$ with
$K_N\subset K_{N+1}$ and $\bigcup_NK_N=X_0$.  For each fixed $N$, local
smooth convergence on $K_N$ and nonnegativity of the regularized density give
\[
 \int_{K_N}
 \bigl|\nabla^{\mathrm{HS}}(t)\Phi(t)\bigr|^2\ddV
 \leq\liminf_{j\to\infty}\cI_{\epsilon_j}(t).
\]
Monotone convergence as $N\to\infty$ then gives
\[
 \cI(t)\leq\liminf_{j\to\infty}\cI_{\epsilon_j}(t).
\]
For $t\geq\tau+1$, integrate the estimate \eqref{eq:one-sided-growth}
from Lemma~\ref{lem:regularized-dissipation-growth} over
$s\in[t-1,t]$, take the lower limit, and use
\eqref{eq:window-convergence}.  This proves
\eqref{eq:dissipation-decay}; its right side tends to zero by
\eqref{eq:total-dissipation}.
\end{proof}

\section{Subsequential convergence of the Y--M--H flow}
\label{sec:slice-compactness}

\subsection{The gauge formulation of the flow}

We use the initial metric $h_0=h(0)$ defined in
Section~\ref{sec:heat-kernel-estimates} as the fixed reference metric.
The metric $h(t)$ is smooth on $X_0$ up to $t=0$,
so on every compact time interval the fiberwise ODE
\begin{equation}\label{eq:horizontal-gauge}
 \dot g(t)=-g(t)\Phi(t),\quad g(0)=\Id
\end{equation}
has a unique smooth complex-gauge solution.  If
$b=h_0^{-1}h$, then
\[
 \dot b=-2b\Phi,\quad \Phi^{*h_0}b=b\Phi.
\]
The endomorphism $c=g^{*h_0}g$ satisfies
\[
 \dot c=-\Phi^{*h_0}c-c\Phi,\quad c(0)=\Id,
\]
and $b$ satisfies the same linear ODE with the same initial value.
Consequently,
\begin{equation}\label{eq:horizontal-isometry}
 h(t)(v,w)=h_0(g(t)v,g(t)w).
\end{equation}

Define
\[
 \nabla_{A(t)}=g(t)\nabla_{h(t)}g(t)^{-1},\quad
 \theta(t)=g(t)\theta g(t)^{-1},\quad
 \Psi(t)=g(t)\Phi(t)g(t)^{-1}.
\]
We retain $\theta$ for the original, time-independent Higgs field and write
$\theta(t)$ for its transport to the fixed metric $h_0$, as in the
Introduction; in particular, $\theta(0)=\theta$.
Let $\nabla_{A(t)}^*$ denote the formal adjoint of the exterior
covariant derivative $\nabla_{A(t)}$ with respect to the $L^2$ inner
product induced by $h_0$ and $\omega$.
Then $A(t)$ is unitary for the fixed metric $h_0$ and $\Psi(t)$ is
$h_0$-self-adjoint.  A unitary connection--Higgs pair $(A,\varphi)$
consists of a connection $A$ on $\cE|_{X_0}$ preserving $h_0$ and a form
$\varphi\in\Omega^{1,0}(X_0,\operatorname{End}(\cE|_{X_0}))$ satisfying
\[
 \Theta_A^{0,2}=0,\quad
 \bar\partial_A\varphi=0,\quad
 \varphi\wedge\varphi=0,
\]
where $\Theta_A$ denotes the curvature of $A$.
For such a pair, set
\[
 \Theta^{1,1}_{A,\varphi}
 =\Theta_A+[\varphi,\varphi^{*h_0}].
\]
Gauge transport gives
\begin{gather}
 \Psi=\ii\Lambda_\omega\Theta^{1,1}_{A,\theta(t)}-\lambda\Id,\quad
 \Theta_A=g\Theta_hg^{-1},\quad
 \theta(t)^{*h_0}=g\theta^{*h}g^{-1},\label{eq:gauge-tensor-transport}\\
 \partial_A\theta(t)=g(\partial_h\theta)g^{-1},\quad
 \nabla_A\Psi=g(\nabla_h\Phi)g^{-1},\quad
 [\Psi,\theta(t)]=g[\Phi,\theta]g^{-1}.\nonumber
\end{gather}
In particular, corresponding tensors in these identities have equal
pointwise norms in the fixed and moving metrics.

Differentiating the transformed connection gives
\[
 \frac{d}{dt}(g\nabla_hg^{-1})
 =g\dot\nabla_hg^{-1}
  -\nabla_A(\dot gg^{-1}).
\]
The metric equation gives
$\dot\nabla_h=\partial_h(h^{-1}\dot h)=-2\partial_h\Phi$, while
$\dot gg^{-1}=-\Psi$.  Hence
\begin{equation}\label{eq:unitary-YMH-flow}
 \dot A=-(\partial_A-\bar\partial_A)\Psi,\quad
 \dot\theta(t)=-[\Psi,\theta(t)].
\end{equation}
The Kähler identities
\[
 \partial_A^*=\ii[\Lambda_\omega,\bar\partial_A],\quad
 \bar\partial_A^*=-\ii[\Lambda_\omega,\partial_A]
\]
and the Bianchi identity give
$\nabla_A^*\Theta_A=\ii(\partial_A-\bar\partial_A)
\Lambda_\omega\Theta_A$.  Thus \eqref{eq:unitary-YMH-flow} is equivalently
\[
 \begin{split}
 \dot A={}&-\nabla_A^*\Theta_A
 -\ii(\partial_A\Lambda_\omega-\bar\partial_A\Lambda_\omega)
       [\theta(t),\theta(t)^{*h_0}],\\
 \dot\theta(t)={}&-
 \bigl[\ii\Lambda_\omega
       (\Theta_A+[\theta(t),\theta(t)^{*h_0}]),\theta(t)\bigr].
 \end{split}
\]
The isometry \eqref{eq:horizontal-isometry} and
\eqref{eq:gauge-tensor-transport} yield
\begin{align}
 \cI(t)
 &=\int_{X_0}\bigl(|\nabla_A\Psi|^2+2|[\Psi,\theta(t)]|^2\bigr)\ddV,\label{eq:unitary-dissipation}\\
 |\dot A|^2+|\dot\theta(t)|^2+|\dot\theta(t)^{*h_0}|^2
 &=|\nabla_A\Psi|^2+2|[\Psi,\theta(t)]|^2.\label{eq:unitary-velocity}
\end{align}
The second identity uses self-adjointness of $\Psi$, which makes the two
Chern-type components of $\nabla_A\Psi$ adjoints of equal norm.

Each regularized flow admits the analogous horizontal transport
\[
 \dot g_\epsilon=-g_\epsilon\Phi_\epsilon,\quad
 g_\epsilon(0)=\Id,
\]
with reference metric $\widehat h$.  The prescribed local
$C^\infty$ convergence $h_{\epsilon_j}\to h$, including at time zero,
and smooth dependence of solutions on the ODE coefficients imply
$g_{\epsilon_j}\to g$, together with convergence of all transformed
tensors, projections, and quotient isometries on every compact finite-time
cylinder. 

\subsection{Construction of the limiting Higgs bundle}

Before constructing the limiting Higgs bundle, we establish the following
estimate using small-energy regularity.

Fix a prescribed time sequence $t_i\to\infty$ and write
\begin{equation}\label{eq:full-density}
 e_2(t)=|\Theta_{A(t)}|^2+2|\partial_{A(t)}\theta(t)|^2.
\end{equation}
The three bidegree components of the full Hitchin--Simpson curvature are
$\Theta^{1,1}_{h,\theta}=\Theta_h+[\theta,\theta^{*h}]$,
$\partial_h\theta$, and
$\bar\partial_{\cE}(\theta^{*h})$.  They are pairwise orthogonal, and the
last two are Hermitian adjoints.  Hence
\begin{equation}\label{eq:full-curvature-decomposition}
 |\Theta^{\mathrm{HS}}_{h,\theta}|^2
 =|\Theta^{1,1}_{h,\theta}|^2+2|\partial_h\theta|^2.
\end{equation}
By \eqref{eq:gauge-tensor-transport}, the corresponding norms satisfy
$|\partial_h\theta|_{h(t),\omega}=|\partial_A\theta(t)|_{h_0,\omega}$.
Proposition~\ref{prop:positive-flow} and
\eqref{eq:full-curvature-decomposition} give uniform positive-time
$L^2$ bounds for
$\Theta_{A(t)}+[\theta(t),\theta(t)^{*h_0}]$
and $\partial_{A(t)}\theta(t)$.
The bound for $\Theta_{A(t)}$ follows from
\eqref{eq:positive-time-bounds} and \eqref{eq:gauge-tensor-transport}.
Consequently, for every $\tau>0$,
\begin{equation}\label{eq:uniform-full-energy}
 \sup_{t\geq\tau}\int_{X_0}e_2(t)\ddV<\infty.
\end{equation}

We use the tensor product of the connection induced by $A(t)$ on
$\operatorname{End}(\cE|_{X_0})$ and the Chern connection induced by
$\omega$ on $T^{*1,0}X$ to differentiate $\theta(t)$.
We denote this covariant derivative by
$\nabla_{A,T^{*1,0}X}\theta(t)$. For iterated derivatives, we write
$\nabla_A^k$, including the fixed connections on all cotangent factors.
The norms in the following lemma are induced by $h_0$ and $\omega$. Set
\[
 P_a(x,t)=B_a(x)\times[t-a^2,t+a^2].
\]

\begin{lemma}[small energy regularity]\label{lem:slice-cylinder}
For every $\tau>0$, there are constants $R_0,\kappa,\delta_0>0$,
depending only on $\tau$ and the fixed geometry, regularization, and initial
data, with the following property.
Suppose that $B_{2R}(x)\Subset X_0$, $0<R\leq R_0$, and
\[
 \limsup_{i\to\infty}R^{4-2n}
 \int_{B_R(x)}e_2(t_i)\ddV<\kappa.
\]
Set $r=R/8$. For every integer $k\geq0$ and every fixed
$0<\delta<\delta_0$, we have
\begin{equation}\label{eq:slice-higher-derivatives}
 \sup_{P_{\delta r}(x,t_i)}
 \left(
 |\nabla_A^k\Theta_A|+
 |\nabla_A^{k+1}\theta|
 \right)\leq C_{k,x,\tau,R,\delta}
\end{equation}
for all sufficiently large $i$, where the constant is independent of $i$.
\end{lemma}

\begin{proof}
The proof follows the argument of
\cite[proof of Theorem~3.2]{LZ11}, with minor modifications using
\cite[Theorem~2.6]{LZ17} and
\cite[Propositions~2.2 and~2.4]{LZZ17Survey}.
We omit the details.
\end{proof}

For the remainder of this subsection, fix $\tau=1$ and choose
$R_0,\kappa,\delta_0$ as in Lemma~\ref{lem:slice-cylinder}.
Let $\mu_i=e_2(t_i)dV_{\omega}$ on $X_0$, extended by zero across
$\Sigma_{\mathcal E}$.  Their total masses are uniformly bounded by
\eqref{eq:uniform-full-energy}.  We first treat $n\geq2$.
Choose a countable dense set
$\mathcal D\subset X$ and consider all balls $B_q(p)$ with
$p\in\mathcal D$ and positive rational $q$.  A diagonal extraction gives
one subsequence for which $\mu_i(B_q(p))$ converges
for every such ball.  If $x\in X$ and $0<r<R$, choose an intermediate
rational ball
\[
 B_r(x)\subset B_q(p)\subset B_R(x).
\]
The convergence of the mass of this intermediate ball gives
\begin{equation}\label{eq:ball-mass-sandwich}
 \limsup_{i\to\infty}\mu_i(B_r(x))
 \leq\liminf_{i\to\infty}\mu_i(B_R(x)).
\end{equation}

Put $\ell=2n-4$ and define
\begin{equation}\label{eq:concentration-set}
 \Sigma_{\mathrm{an}}
 =\left\{x\in X_0:
 \liminf_{r\downarrow0}r^{-\ell}
 \limsup_{i\to\infty}\mu_i(B_r(x))\geq\kappa\right\}.
\end{equation}
If $x\notin\Sigma_{\mathrm{an}}$, there are
$\eta_x<\kappa$ and radii $r_j\downarrow0$ such that
\[
 r_j^{-\ell}\limsup_{i\to\infty}\mu_i(B_{r_j}(x))
 \leq\eta_x.
\]
If $x\in\Sigma_{\mathrm{an}}$, then for all sufficiently small $r$,
the definition at radius $r/2$ and
\eqref{eq:ball-mass-sandwich} give
\[
 \liminf_{i\to\infty}\mu_i(B_r(x))
 \geq\kappa\,2^{-\ell-1}r^\ell.
\]
Consequently, for every finite pairwise disjoint family of such balls,
\[
 \kappa\,2^{-\ell-1}\sum_a r_a^\ell
 \leq\sum_a\liminf_i\mu_i(B_{r_a}(x_a))
 \leq\liminf_i\sum_a\mu_i(B_{r_a}(x_a))
 \leq\sup_i\mu_i(X).
\]
The $5r$ selection lemma applied to the resulting fine cover therefore
gives
$\mathcal H^{2n-4}(\Sigma_{\mathrm{an}})<\infty$.  When $n=2$, the same
finite-family estimate bounds the number of selected balls, and hence shows
that $\Sigma_{\mathrm{an}}$ is finite.

Fix $x\in X_0\setminus\Sigma_{\mathrm{an}}$.  Choose one of the preceding
radii $R\leq R_0$ such that $B_{2R}(x)\Subset X_0$.
Since $\omega$ is Kähler,
$\partial_{A(t)}\theta(t)$ is obtained by antisymmetrizing the $(1,0)$
part of $\nabla_{A,T^{*1,0}X}\theta(t)$, so
\[
 |\partial_{A(t)}\theta(t)|
 \leq C\,|\nabla_{A(t),T^{*1,0}X}\theta(t)|.
\]
The strict normalized limsup bound allows us to apply
Lemma~\ref{lem:slice-cylinder}.  Together with the definition
\eqref{eq:full-density}, it gives a neighborhood $W\Subset X_0$ of $x$
and a constant $C_W$ such that $e_2(t_i)\leq C_W$ on $W$ for all
sufficiently large $i$.  Thus, for $y$ in a smaller neighborhood and all
sufficiently small $a$,
\[
 a^{4-2n}\limsup_i\mu_i(B_a(y))\leq C_Wa^4\longrightarrow0.
\]
The complement of $\Sigma_{\mathrm{an}}$ is therefore relatively open in
$X_0$.

By reflexivity, $\Sigma_{\mathcal E}$ is analytic of complex codimension
at least three, hence closed and of zero $\mathcal H^{2n-4}$ measure.  Therefore
\begin{equation}\label{eq:bad-set}
 Z=\Sigma_{\mathcal E}\cup\Sigma_{\mathrm{an}}
\end{equation}
is closed and has finite $\mathcal H^{2n-4}$ measure.

When $n=1$, reflexivity gives $\Sigma_{\mathcal E}=\varnothing$, and
the uniform mass bound gives, for every $x\in X$,
\[
 r^2\limsup_i\mu_i(B_r(x))
 \leq r^2\sup_i\mu_i(X)\longrightarrow0\quad(r\downarrow0).
\]
Thus every point satisfies the strict small-energy hypothesis of
Lemma~\ref{lem:slice-cylinder} for all sufficiently small radii.
In this case we set $\Sigma_{\mathrm{an}}=Z=\varnothing$.
The following construction applies in both cases.

\begin{proposition}\label{prop:slice-compactness}
For the prescribed $t_i\to\infty$, there are a subsequence, relabelled
$t_i$, and a Hermitian bundle $(E_{\infty},h_\infty)$ on $U=X\setminus Z$.
Here $Z$ is the closed set constructed above; it contains
$\Sigma_{\mathcal E}$, satisfies $\mathcal H^{2n-4}(Z)<\infty$ when
$n\geq2$, and is empty when $n=1$.
There are a smooth $h_\infty$-unitary connection $A_\infty$ on
$E_\infty$ and a smooth form
$\theta_\infty\in\Omega^{1,0}(U,\operatorname{End}E_\infty)$.
A countable locally finite cover $\{U_\alpha\}$ of $U$ admits
$h_0$-unitary gauges $a_{i,\alpha}$ on $\cE|_{U_\alpha}$, the same
for the connection and the Higgs field, such that
\begin{equation}\label{eq:slice-convergence}
 a_{i,\alpha}\bigl(A(t_i),\theta(t_i)\bigr)
 \longrightarrow(A_{\infty,\alpha},\theta_{\infty,\alpha})
\end{equation}
in $C^\infty$ on compact subsets of $U_\alpha$, where the right-hand side
is the local expression of $(A_\infty,\theta_\infty)$ in an
$h_\infty$-unitary trivialization of $E_\infty|_{U_\alpha}$.
The metric $h_\infty$ is obtained by gluing the standard Hermitian
metrics in these local trivializations using the limiting unitary
transition functions constructed in the proof.  Moreover,
\begin{equation}\label{eq:limiting-Higgs-equations}
 \Theta_{A_\infty}^{0,2}=0,\quad
 \bar\partial_{A_\infty}\theta_\infty=0,\quad
 \theta_\infty\wedge\theta_\infty=0.
\end{equation}
If
\[
 \Psi_\infty=\ii\Lambda_\omega
 \bigl(\Theta_{A_\infty}
       +[\theta_\infty,\theta_\infty^{*h_\infty}]\bigr)-\lambda\Id,
\]
then
\begin{equation}\label{eq:critical-slice}
 \nabla_{A_\infty}\Psi_\infty=0,\quad
 [\Psi_\infty,\theta_\infty]=0.
\end{equation}
The limiting pair satisfies
\[
 \Theta_{A_\infty}\in L^2(U),\quad
 \Lambda_\omega\Theta_{A_\infty}\in L^\infty(U),
\]
with norms taken with respect to $h_\infty$ and $\omega$, and
integration against $dV_\omega$.
\end{proposition}

\begin{proof}
Write $A_i=A(t_i)$, $\theta_i=\theta(t_i)$, and
$r_E=\rank\cE$.
For every $x\in U$, small energy regularity gives a coordinate ball $V_x\Subset U$
on which
\begin{equation}\label{eq:local-covariant-bounds}
 \sup_{V_x}\left(
 |\nabla_{A_i}^k\Theta_{A_i}|+|\nabla_{A_i}^k\theta_i|
 \right)\leq C_{k,x},\quad k\geq0,
\end{equation}
for all sufficiently large $i$.  All derivatives here include the fixed
connections on the tensor factors.  Choose a countable locally finite
cover $\{U_\alpha\}$ of $U$ by coordinate balls, with slightly larger
concentric coordinate balls $V_\alpha$ such that
$\overline{U_\alpha}\subset V_\alpha\Subset U$ and
\eqref{eq:local-covariant-bounds} holds on each $V_\alpha$.

Fix $q\in(n,2n)$. By the uniform curvature bound in
\eqref{eq:local-covariant-bounds}, we may choose the preceding
cover sufficiently fine that, after rescaling each $V_\alpha$
to the unit ball, the curvature is uniformly small in $L^q$
for all sufficiently large $i$. Choose an $h_0$-unitary trivialization on each $V_\alpha$.
Uhlenbeck's local gauge theorem
\cite[Theorem~1.3]{Uhlenbeck82} provides local
$h_0$-unitary gauges $a_{i,\alpha}$ such that for
\[
 d+A_{i,\alpha}
 =a_{i,\alpha}\nabla_{A_i}a_{i,\alpha}^{-1},
 \quad
 \theta_{i,\alpha}
 =a_{i,\alpha}\theta_i a_{i,\alpha}^{-1},
\]
we have
\[
 d^*A_{i,\alpha}=0,
 \quad
 \|A_{i,\alpha}\|_{W^{1,q}(V_\alpha)}
 \leq C_\alpha.
\]
The constants are independent of $i$.

Put
$\Theta_{i,\alpha}
=a_{i,\alpha}\Theta_{A_i}a_{i,\alpha}^{-1}$.
The curvature identity
\[
 dA_{i,\alpha}
 =\Theta_{i,\alpha}
  -A_{i,\alpha}\wedge A_{i,\alpha},
\]
together with the Coulomb condition and
\eqref{eq:local-covariant-bounds}, yields uniform local
Sobolev bounds of all orders for $A_{i,\alpha}$ by interior
elliptic estimates. The corresponding estimates for
$\theta_{i,\alpha}$ follow from its covariant derivative
bounds in the same gauges. Sobolev embedding therefore gives
\begin{equation}\label{eq:local-gauge-bounds}
 \|A_{i,\alpha}\|_{C^k(U_\alpha)}
 +\|\theta_{i,\alpha}\|_{C^k(U_\alpha)}
 \leq C_{k,\alpha},
 \quad k\geq0,
\end{equation}
with constants independent of $i$. The Arzel\`a--Ascoli theorem and the diagonal method provide a subsequence
\begin{equation}\label{eq:local-pair-limit}
 (A_{i,\alpha},\theta_{i,\alpha})
 \longrightarrow(A_{\infty,\alpha},\theta_{\infty,\alpha})
 \quad\text{in }C^\infty_{\mathrm{loc}}(U_\alpha).
\end{equation}

To construct a bundle from these local limits, we also take limits of
the transition functions.  Let
$g_{\alpha\beta}:U_\alpha\cap U_\beta\to\mathrm U(r_E)$ be the
transition functions of the fixed trivializations, with fiber coordinates
related by $v_\alpha=g_{\alpha\beta}v_\beta$.  In the frames just
constructed, the transition functions are
\begin{equation}\label{eq:gauged-transition-functions}
 g_{i,\alpha\beta}
 =a_{i,\alpha}g_{\alpha\beta}a_{i,\beta}^{-1}.
\end{equation}
Since the local connection forms and Higgs fields represent the same
objects on $\cE|_U$, they satisfy
\begin{equation}\label{eq:transition-compatibility}
 \begin{split}
 dg_{i,\alpha\beta}
 &=g_{i,\alpha\beta}A_{i,\beta}
   -A_{i,\alpha}g_{i,\alpha\beta},\\
 \theta_{i,\alpha}g_{i,\alpha\beta}
 &=g_{i,\alpha\beta}\theta_{i,\beta}.
 \end{split}
\end{equation}
Unitarity bounds $g_{i,\alpha\beta}$ pointwise.  Differentiating the
first equation and using \eqref{eq:local-gauge-bounds} bounds all of
its derivatives on compact subsets of $U_\alpha\cap U_\beta$.
A further diagonal subsequence over these overlaps and their compact
exhaustions therefore gives
\[
 g_{i,\alpha\beta}\longrightarrow g_{\infty,\alpha\beta}
 \quad\text{in }C^\infty_{\mathrm{loc}}(U_\alpha\cap U_\beta).
\]
The limits are unitary and retain the transition identities:
\[
 g_{\infty,\alpha\alpha}=\Id,\quad
 g_{\infty,\beta\alpha}=g_{\infty,\alpha\beta}^{-1},\quad
 g_{\infty,\alpha\beta}g_{\infty,\beta\gamma}
 =g_{\infty,\alpha\gamma}.
\]
The last equality holds on every triple overlap.
Define
\begin{equation}\label{eq:limit-bundle-gluing}
 E_\infty
 =\left(\coprod_\alpha U_\alpha\times\mathbb C^{r_E}\right)\big/\sim,
 \quad
 (x,v)_\beta\sim(x,g_{\infty,\alpha\beta}(x)v)_\alpha.
\end{equation}
The transition identities make $E_\infty$ a smooth vector bundle.
Since its transition functions are unitary, the standard Hermitian
metrics on the local trivial bundles define a metric $h_\infty$ on
$E_\infty$.

Taking limits in \eqref{eq:transition-compatibility} gives
\[
 \begin{split}
 A_{\infty,\alpha}
 &=g_{\infty,\alpha\beta}A_{\infty,\beta}
       g_{\infty,\alpha\beta}^{-1}
   -(dg_{\infty,\alpha\beta})g_{\infty,\alpha\beta}^{-1},\\
 \theta_{\infty,\alpha}
 &=g_{\infty,\alpha\beta}\theta_{\infty,\beta}
       g_{\infty,\alpha\beta}^{-1}.
 \end{split}
\]
Consequently, the local connections $d+A_{\infty,\alpha}$ define an
$h_\infty$-unitary connection $A_\infty$ on $E_\infty$, and the forms
$\theta_{\infty,\alpha}$ define
$\theta_\infty\in\Omega^{1,0}(U,\operatorname{End}E_\infty)$.
Equation~\eqref{eq:local-pair-limit} is precisely the convergence
asserted in \eqref{eq:slice-convergence}.

Every pair $(A_i,\theta_i)$ satisfies
\[
 \Theta_{A_i}^{0,2}=0,\quad
 \bar\partial_{A_i}\theta_i=0,\quad
 \theta_i\wedge\theta_i=0.
\]
These equations are invariant under the unitary changes of frame, so
local smooth convergence gives \eqref{eq:limiting-Higgs-equations}.
The first equation gives $\bar\partial_{A_\infty}^2=0$, so
$\bar\partial_{A_\infty}$ defines a holomorphic structure on
$E_\infty$.  The other two say that
$\theta_\infty$ is holomorphic and satisfies the Higgs integrability
condition.  This constructs the Higgs bundle
$(E_\infty,\bar\partial_{A_\infty},\theta_\infty)$ on $U$.

It remains to prove the critical-point equations and the curvature
bounds.  Set $\Psi_i=\Psi(t_i)$.  Its expression in each of the same
unitary frames converges smoothly to the local expression of
$\Psi_\infty$, by \eqref{eq:gauge-tensor-transport} and
\eqref{eq:local-pair-limit}.  The squared norms in
\eqref{eq:unitary-dissipation} are invariant under these changes of
frame and thus converge uniformly on every compact $K\Subset U$.
Proposition~\ref{prop:dissipation} gives
\[
 \int_K\left(
 |\nabla_{A_\infty}\Psi_\infty|^2
 +2|[\Psi_\infty,\theta_\infty]|^2\right)\ddV
 \leq\lim_{i\to\infty}\cI(t_i)=0,
\]
which proves \eqref{eq:critical-slice}.  Local smooth convergence,
\eqref{eq:uniform-full-energy}, Fatou's lemma, and exhaustion of $U$ give
finite $L^2$ energy for the Chern curvature.  The positive-time bounds for $\Phi(t)$ and $\theta$ pass locally to the limit, and
\[
 \Lambda_\omega\Theta_{A_\infty}
 =-\ii(\Psi_\infty+\lambda\Id)
  -\Lambda_\omega[\theta_\infty,\theta_\infty^*]
\]
gives the asserted bound for the contracted Chern curvature.
\end{proof}

\section{Reflexive extension and spectral decomposition of the limit Higgs sheaf}
\label{sec:extension-polystability}

\begin{proposition}\label{prop:extension}
The Higgs bundle
$(E_{\infty},\bar\partial_{A_\infty},\theta_\infty)$ on $U$ has a unique
reflexive Higgs extension $(\cE_{\infty},\theta_{\infty})$ on $X$.
The eigenvalues of $\Psi_\infty$ and their multiplicities are constant on
$U$, and their eigensubbundles extend to a decomposition
\begin{equation}\label{eq:spectral-decomposition}
 (\cE_{\infty},\theta_{\infty})
 =\bigoplus_{c\in\operatorname{Spec}(\Psi_\infty)}
   (\cE_{\infty,c},\theta_{\infty,c}).
\end{equation}
Here $\cE_{\infty,c}|_U=\ker(\Psi_\infty-c\Id)$.
Each summand is a reflexive Higgs-polystable sheaf and carries a smooth admissible
metric $h_c$ on its locally free locus satisfying
\begin{equation}\label{eq:block-HYM}
 \ii\Lambda_{\omega}
 \Theta^{1,1}_{h_c,\theta_{\infty,c}}
 =(\lambda+c)\Id.
\end{equation}
In particular,
\begin{equation}\label{eq:block-slope}
 \mu_{\omega}(\cE_{\infty,c})
 =\frac{\Vol(X,\omega)}{2\pi}(\lambda+c).
\end{equation}
\end{proposition}

\begin{proof}
\emph{Reflexive extension.}
For $n\geq2$, the Chern curvature bounds and
$\mathcal H^{2n-4}(Z)<\infty$ from
Proposition~\ref{prop:slice-compactness} allow us to apply Bando--Siu's ordinary
extension theorem \cite[Theorem~2(a)]{BS94}.  This gives a reflexive extension
$\cE_\infty$ of the underlying holomorphic bundle.
When $n=1$, the limiting Higgs bundle is already defined on all of $X$,
so no extension is needed.

The Hausdorff measure bound implies that $Z$ has locally zero
$\mathcal H^{2n-2}$ measure.  On the common locally free locus of two
reflexive sheaves, Shiffman's theorem
\cite[Lemma~3, p.~115]{Shi68} extends the coefficients of a holomorphic
bundle map uniquely across $Z$.  Across the remaining analytic singular
sets, the map extends as a section of a reflexive Hom sheaf
\cite[Propositions~5.5.21 and~5.5.23]{Kobayashi87}.  Transition,
intertwining, integrability, and inverse identities extend because both
sides are sections of torsion-free Hom sheaves and agree on a dense open set.  Thus a
holomorphic Higgs map defined on $U$ extends uniquely between any prescribed
reflexive extensions, and an isomorphism extends together with its inverse.

Apply this principle to
$\theta_\infty:E_\infty\to E_\infty\otimes\Omega_X^1$.  It extends uniquely
to a global Higgs field on $\cE_\infty$; the extension is
integrable because $\theta_\infty\wedge\theta_\infty=0$ on $U$.  The same
principle proves uniqueness of the resulting reflexive Higgs extension
relative to its prescribed identification with $(E_\infty,\theta_\infty)$
on $U$.

The metric conclusion of
\cite[Theorem~2(b)]{BS94}, applied separately to the limiting bundle and
its Hermitian dual, bounds the relative metric and its inverse locally on
$X\setminus\Sing(\cE_\infty)$.  The continuation principle identifies the
dual extension with $\cE_\infty^\vee$.  Hence the limiting metric is locally
uniformly comparable in both directions with every smooth reference metric
on the locally free locus.

\smallskip
\noindent\emph{Spectral decomposition and regularity.}
For $n\geq2$, $\mathcal H^{2n-4}(Z)<\infty$ implies
$\mathcal H^{2n-1}(Z)=0$. It's easy to see that the complement $U$ is connected.

Fix $x_0\in U$.  Parallel transport along a path $\gamma$ from $x_0$ to
$x$ satisfies
\[
 \Psi_\infty(x)T_\gamma=T_\gamma\Psi_\infty(x_0)
\]
because $\nabla_{A_\infty}\Psi_\infty=0$.  Thus the real eigenvalues and
their multiplicities are constant on $U$.  For each
$c\in\operatorname{Spec}(\Psi_\infty)$, the eigensubbundle
$E_{\infty,c}=\ker(\Psi_\infty-c\Id)$ is smooth and $A_\infty$-parallel.
These subbundles are orthogonal because $\Psi_\infty$ is self-adjoint,
holomorphic because $A_\infty^{0,1}$ preserves them, and Higgs invariant
because $[\Psi_\infty,\theta_\infty]=0$.  They give an orthogonal parallel
holomorphic Higgs splitting, and restriction of the definition of
$\Psi_\infty$ gives \eqref{eq:block-HYM}.

Since the splitting is parallel,
\[
 |\Theta_{A_{\infty,c}}|\leq|\Theta_{A_\infty}|,\quad
 |\Lambda_\omega\Theta_{A_{\infty,c}}|
 \leq|\Lambda_\omega\Theta_{A_\infty}|,\quad
 |\theta_{\infty,c}|\leq|\theta_\infty|.
\]
Bando--Siu's theorem \cite[Theorem~2(a)]{BS94} therefore extends every
block to a reflexive sheaf $\cE_{\infty,c}$.  The continuation principle
extends its Higgs field uniquely.  The finite direct sum of these extensions
agrees with $(E_\infty,\theta_\infty)$ on $U$, so uniqueness of the extension
identifies this sum with $(\cE_\infty,\theta_\infty)$ and proves
\eqref{eq:spectral-decomposition}.

Let $Y_c=X\setminus\Sing(\cE_{\infty,c})$ and choose a smooth reference
metric $h_{c,0}$ there.  For
\[
 H_c=h_{c,0}^{-1}h_c,
\]
\cite[Theorem~2(b)]{BS94} gives
$H_c\in W^{2,p}_{\mathrm{loc}}(Y_c)$ for every finite $p$ and local
boundedness of $H_c$.  Applying the theorem separately to the dual block
and using the continuation principle to identify its extension with
$\cE_{\infty,c}^{\vee}$ gives local boundedness of $H_c^{-1}$.

The set $Z$ has zero volume, so the block equation holds weakly on $Y_c$:
\[
 \ii\Lambda_\omega\left(
 \Theta_{h_{c,0}}
 +\bar\partial\bigl(H_c^{-1}\partial_{h_{c,0}}H_c\bigr)
 +[\theta_{\infty,c},
   H_c^{-1}\theta_{\infty,c}^{*h_{c,0}}H_c]\right)
 =(\lambda+c)\Id.
\]
Here
$\Theta_{h_c}=\Theta_{h_{c,0}}+
\bar\partial(H_c^{-1}\partial_{h_{c,0}}H_c)$ and
$\theta_{\infty,c}^{*h_c}
=H_c^{-1}\theta_{\infty,c}^{*h_{c,0}}H_c$.
Choose $p>2n$.  Sobolev embedding gives a
$C^{1,\alpha}_{\mathrm{loc}}$ representative of $H_c$, while the inverse
bound gives a positive lower bound for the eigenvalues, which persists by
continuity from $Y_c\cap U$ to all of $Y_c$.  After multiplication by $H_c$, the weak equation has a fixed componentwise
scalar elliptic operator as its principal part; every remaining term contains only smooth reference coefficients,
$H_c$, $H_c^{-1}$, and at most first derivatives of $H_c$.  Schauder
estimates first give $C^{2,\alpha}_{\mathrm{loc}}$ regularity and then,
by iteration, smoothness.  See also
\cite[equations~(2.37)--(2.42) and the paragraph preceding
Proposition~2.10]{LZ17}.  Adding $Y_c\cap Z$ does not change the integral of
$|\Theta_{h_c}|^2$, and continuity preserves the bound for
$|\Lambda_\omega\Theta_{h_c}|$.  Thus $h_c$ is
smooth and admissible on $Y_c$, and \eqref{eq:block-HYM} holds there
pointwise.

Taking the trace of (6.2) and applying Chern--Weil theory to the determinant metric, which
extends smoothly by \cite[Theorem 2(c)]{BS94}, gives (6.3). Since $h_c$ is admissible and satisfies (6.2), the Higgs-polystability of $(\cE_{\infty,c},\theta_{\infty,c})$ can be obtained by standard arguments (see for example \cite[Corollary 3.5]{BiswasSchumacher09}).
\end{proof}

\section{Harder--Narasimhan type and the limiting spectral vector}
\label{sec:hn-vector}

In this section, we compare the ranks and slopes of summands of $\Gr^{\HNS}(\mathcal E,\theta)^{\vee\vee}$ with those of $(\cE_{\infty},\theta_{\infty})=\bigoplus_c(\cE_{\infty,c},\theta_{\infty,c})$. Let $r_s$ be the rank of
$\cE_s/\cE_{s-1}$, and set
\begin{equation*}
 R=\rank\cE=\sum_{s=1}^{m}r_s,\quad
 R_s=\sum_{q=1}^{s}r_q,\quad R_0=0.
\end{equation*}
Put
\begin{equation}\label{eq:HNS-vector}
 a_s=\frac{2\pi}{\Vol(X,\omega)}\mu_{\omega}(\cE_s/\cE_{s-1}),
\end{equation}
with $\boldsymbol{\mu}=(\mu_1,\ldots,\mu_R)$ the nonincreasing vector obtained by repeating
$a_s$ exactly $r_s$ times.  Let $b=(b_1,\ldots,b_R)$ be the nonincreasing vector obtained
by repeating $\lambda+c$ exactly $\mathrm{rank} (\cE_{\infty,c})$ times.

\begin{lemma}\label{lem:polygon}
For $M_k=\sum_{q=1}^{k}\mu_q$ and $B_k=\sum_{q=1}^{k}b_q$,
\begin{equation}\label{eq:majorization}
 M_k\leq B_k\quad(1\leq k<R),\quad M_R=B_R.
\end{equation}
\end{lemma}

\begin{proof}
Let $\mathcal F\subset\cE$ be a saturated $\theta$-invariant subsheaf of rank
$q$. On $Y$, put
\[
 \widetilde{\mathcal F}
 =\Sat_{\cB}\im(\rho^*\mathcal F\longrightarrow \rho^*\mathcal E \longrightarrow\cB).
\]
The image is $\beta$-invariant, and so is its saturation $\widetilde{\mathcal F}$.

Firstly we show the convergence of $\deg_{\omega_\epsilon}(\widetilde{\mathcal F})$ as $\epsilon \to 0$. For the three determinant line bundles $\det \widetilde{\mathcal F}$, $\det(\rho^* \mathcal F)$, and $\rho^*(\det \mathcal F)$ defined on $Y$, it is easy to see that both $\det \widetilde{\mathcal F} \otimes (\det(\rho^* \mathcal F))^{-1}$ and $\det(\rho^* \mathcal F) \otimes (\rho^*(\det \mathcal F))^{-1}$ are tensor products of line bundles, each being either the determinant line bundle of a torsion sheaf or its dual, and so is $\det \widetilde{\mathcal F} \otimes (\rho^*(\det \mathcal F))^{-1}$. It is also known that these torsion sheaves vanish on the biholomorphic locus $Y \setminus D$. By \cite[Proposition 5.6.14]{Kobayashi87}, we may construct a global meromorphic section $s_{\mathcal F}$ on $\det \widetilde{\mathcal F} \otimes (\rho^*(\det \mathcal F))^{-1}$ whose divisor $D_{\mathcal F}$ is only supported in $D$. The Lelong--Poincaré equation gives
\begin{equation*}
 c_1(\widetilde{\mathcal F})=\rho^*c_1(\mathcal F)+[D_{\mathcal F}].
\end{equation*}
Since $\rho(\operatorname{Supp}(D_{\mathcal F}))\subset \Sigma_{\mathcal E}$, which has complex codimension at least three, for any irreducible component $A$ of $D_{\mathcal F}$,
\[
 \int_Y[A]\wedge(\rho^*\omega)^{n-1}:=\int_{A_{\mathrm{reg}}}(\rho^*\omega)^{n-1}=0 \Rightarrow \int_Y[D_{\mathcal F}]\wedge(\rho^*\omega)^{n-1}=0.
\]
Expanding $(\rho^*\omega+\epsilon\eta)^{n-1}$ as a polynomial in $\epsilon$, we have
\[
 \begin{split}
 \deg_{\omega_\epsilon}(\widetilde{\mathcal F})
 &=\int_Y\bigl(\rho^*c_1(\mathcal F)+[D_{\mathcal F}]\bigr)
       \wedge\frac{(\rho^*\omega+\epsilon\eta)^{n-1}}{(n-1)!}\\
 &=\deg_\omega(\mathcal F)+O(\epsilon).
 \end{split}
\]
Consequently, $\deg_{\omega_\epsilon}(\widetilde{\mathcal F})
 \longrightarrow\deg_\omega(\mathcal F)$.

Fix $\epsilon>0$ and $t>0$, and put
\[
 K_\epsilon(t)=\ii\Lambda_{\omega_\epsilon}
 \Theta^{1,1}_{h_\epsilon(t),\beta}.
\]
The endomorphism $K_\epsilon(t)$ is self-adjoint with respect to
$h_\epsilon(t)$.  Let
$\kappa_1(K_\epsilon(t))\geq\cdots\geq\kappa_R(K_\epsilon(t))$
be its ordered real eigenvalues, and let
$P_{\widetilde{\mathcal F},h_\epsilon(t)}$ be the weak
$h_\epsilon(t)$-orthogonal projection onto $\widetilde{\mathcal F}$. The Chern--Weil identity \cite[Lemma~3.2]{Simpson88} gives
\[
 \begin{split}
 2\pi\deg_{\omega_\epsilon}(\widetilde{\mathcal F})
 &=\int_Y\tr\bigl(K_\epsilon(t)
          P_{\widetilde{\mathcal F},h_\epsilon(t)}\bigr)
          dV_{\omega_\epsilon}\\
 &\quad-\int_Y\left(
 |\bar\partial P_{\widetilde{\mathcal F},h_\epsilon(t)}|^2
 +|[\beta,P_{\widetilde{\mathcal F},h_\epsilon(t)}]|^2\right)
 dV_{\omega_\epsilon}.
 \end{split}
\]
A purely linear-algebraic fact yields $\tr\bigl(K_\epsilon(t) P_{\widetilde{\mathcal F},h_\epsilon(t)}\bigr) \leq \sum_{v=1}^{q}\kappa_v(K_\epsilon(t))$. Thus
\begin{equation}\label{eq:regularized-filtration-ky-fan}
 2\pi\deg_{\omega_\epsilon}(\widetilde{\mathcal F})
 \leq\int_Y\sum_{v=1}^{q}\kappa_v(K_\epsilon(t))
 dV_{\omega_\epsilon}.
\end{equation}
For $\mathcal F=\cE$, the saturation is $\cB$, the projection is the identity, and
the derivative term vanishes. In this case
\eqref{eq:regularized-filtration-ky-fan} is an equality.

Keep $t>0$ fixed and set
$K(t)=\ii\Lambda_\omega\Theta^{1,1}_{h(t),\theta}$, which is self-adjoint
with respect to $h(t)$.  Let $\kappa_v(K(t))$ be its ordered eigenvalues.
Let $\epsilon=\epsilon_j\downarrow0$ along the sequence defining the
flow $h(t)$.  The local smooth convergence on $Y_0$
\[
 K_{\epsilon_j}(t)  \longrightarrow  \widetilde K(t):=\ii\Lambda_{\rho^* \omega}\Theta^{1,1}_{\widetilde h(t),\beta},\quad
\widetilde h(t)=(\rho|_{Y_0})^*h(t),
\quad
\widetilde K(t)=(\rho|_{Y_0})^*K(t)
\]
gives convergence of the ordered eigenvalues: $\kappa_v(K_{\epsilon_j}(t))  \rightarrow  \kappa_v(\widetilde K(t))$. Moreover, since $\Phi_{\epsilon_j}$ is uniformly bounded and $\lambda_{\epsilon_j}\to\lambda$, all these eigenvalues are also uniformly bounded. In addition $\omega_{\epsilon_j}^n\leq C\eta^n$. Thus dominated convergence gives
\begin{equation}\label{eq:physical-filtration-ky-fan}
 2\pi\deg_\omega(\mathcal F)
 \leq\int_{X_0}\sum_{v=1}^{q}\kappa_v(K(t))\ddV.
\end{equation}
For $\mathcal F=\cE$, (\ref{eq:physical-filtration-ky-fan}) is still an equality.

Let $t_i \to \infty$ (Proposition~\ref{prop:slice-compactness}) in (\ref{eq:physical-filtration-ky-fan}). Since the adjoint action of the unitary gauge $a_{i,\alpha}$ preserves the eigenvalues, the local smooth convergence on $X\setminus Z$ gives the locally uniform convergence $\kappa_v(K(t_i)) \to \kappa_v(K_{\infty})$ with $K_{\infty}=\ii\Lambda_\omega \Theta^{1,1}_{A_\infty,\theta_\infty}$, which is again uniformly bounded since $\Phi(t_i)$ is uniformly bounded by Proposition~\ref{prop:positive-flow}. Thus dominated convergence gives
\[
2\pi\deg_{\omega}(\mathcal F) \leq \int_{X_0}\sum_{v=1}^{q}\kappa_v(K_{\infty})\ddV.
\]
Moreover, Proposition~\ref{prop:extension} shows that the eigenvalues of $K_\infty$ are constant. Taking $\mathcal F=\cE_s$, we obtain
\[
 2\pi\deg_\omega(\cE_s)
 \leq\Vol(X,\omega)\sum_{v=1}^{R_s}b_v
 =\Vol(X,\omega)B_{R_s}.
\]
By definition, the left side equals  $2\pi \sum_{v=1}^s r_v \cdot \mu_{\omega}(\mathcal E_v / \mathcal E_{v-1})=\Vol(X,\omega)M_{R_s}$, so
$M_{R_s}\leq B_{R_s}$.  For $\mathcal F=\cE$ we have the equality $M_R=B_R$.

Finally, $B_k$ is concave because $b$ is nonincreasing, while $M_k$ is linear on each interval $[R_{s-1}, R_s]$. Hence the endpoint inequalities force $M_k\leq B_k$ for all $k$.
\end{proof}

Next we show the existence of an $L^1$-$\delta$-approximate critical Hermitian metric on $\mathcal B$. The proof is similar to that of Proposition~4.5 in \cite{LZZ18}, which treats the $L^p$-case without Higgs fields.

\begin{lemma}
\label{lem:approximate-critical}
For every $\delta>0$, there is a smooth metric $\widehat h_\delta$ on
$\cB$ with the following properties.  Put
$h_\delta=(\rho^{-1})^*(\widehat h_\delta|_{Y_0})$ and
$K_\delta=\ii\Lambda_\omega\Theta^{1,1}_{h_\delta,\theta}$, and let $P_s$
be the almost-everywhere defined $h_\delta$-orthogonal projection onto
$\cE_s$ on $X_f$, with $P_0=0$.  Then
\begin{equation}\label{eq:approximate-critical}
 \left\|K_\delta-\sum_s a_s(P_s-P_{s-1})\right\|
       _{L^1(X_0,\omega,h_\delta)}<\delta.
\end{equation}
Moreover, set $K_{\delta,\epsilon}=\ii\Lambda_{\omega_\epsilon}
\Theta^{1,1}_{\widehat h_\delta,\beta}$ and let $P_{s,\delta}$ be the
almost-everywhere defined $\widehat h_\delta$-orthogonal projection onto
the saturation in $\cB$ of the image of $\rho^*\cE_s$, with
$P_{0,\delta}=0$.  Then
\begin{equation}\label{eq:regularized-residual-convergence}
 \lim_{\epsilon\downarrow0}
 \left\|K_{\delta,\epsilon}
       -\sum_s a_s(P_{s,\delta}-P_{s-1,\delta})\right\|
       _{L^1(Y,\omega_\epsilon,\widehat h_\delta)}
 =
 \left\|K_\delta-\sum_s a_s(P_s-P_{s-1})\right\|
       _{L^1(X_0,\omega,h_\delta)}.
\end{equation}
\end{lemma}

\begin{proof}
Put $\Sigma_f=X\setminus X_f$.  We employ the second regularization $\sigma\colon \widetilde Y\to Y$ so that the quotient Higgs sheaves $(\mathcal R_s,\widetilde\theta_s):=\widetilde{\mathcal E}_s/\widetilde{\mathcal E}_{s-1}$ from the pulled-back filtration are locally free. Set $\pi=\rho\circ\sigma$, $\widetilde\cB=\sigma^*\cB$, and
$\widetilde\beta=\sigma^*\beta$. Fix a Kähler form $\widetilde\eta$ with
$\Vol(\widetilde Y,\widetilde\eta)=1$, and put
\[
 \omega_0=\pi^*\omega,\quad
 \omega_u=\omega_0+u\widetilde\eta,\quad 0<u\leq1.
\]
In particular,
$\int_{\widetilde Y}dV_{\omega_0}=\Vol(X,\omega)$, although $\omega_0$ is degenerate
on the exceptional locus.

Let $Q_s=(\cE_s/\cE_{s-1})^{\vee\vee}$, equipped with the induced Higgs
field. By definition of the HNS filtration, $Q_s$ is $\omega$-stable by \cite{Kobayashi87}. Stability yields a uniform slope gap $d_s$ such that
\[
 \mu_\omega(\mathcal S)
 \leq\mu_\omega(Q_s)-d_s
\]
for every proper saturated Higgs subsheaf $\mathcal S\subset Q_s$.

If $\mathcal A\subset\mathcal R_s$ is a nonzero proper saturated Higgs
subsheaf, we define over $X_f$:
\[
\mathcal S_{\mathcal A}
 :=\Sat_{Q_s}\im(\pi_* \mathcal A \longrightarrow \pi_* \mathcal R_s \longrightarrow Q_s).
\]
Since $\mathcal S_{\mathcal A}=\pi_*\mathcal A$ away from $\Sigma_f$, which has codimension at least $2$, we have $c_1(\mathcal S_{\mathcal A})=c_1(\pi_*\mathcal A)$. Similarly $c_1(Q_s)=c_1(\pi_* \mathcal R_s)$. Define the $\omega_u$-degree on $\widetilde Y$:
\[
 \deg_{\omega_u}(\mathcal F)
 =\int_{\widetilde Y}c_1(\mathcal F)
   \wedge\frac{\omega_u^{n-1}}{(n-1)!}, \quad \mu_{\omega_u}(\mathcal F) = \frac{\mathrm{deg}_{\omega_u}(\mathcal F)}{\mathrm{rank}(\mathcal F)}.
\]
\cite[Theorem~4.9]{Sibley15} gives slope control relating sheaves and their push-downs:
\[
 \mu_{\omega_u}(\mathcal A)
 \leq\mu_\omega(\mathcal S_{\mathcal A})+C_su,\quad
 \mu_{\omega_u}(\mathcal R_s)
 \geq\mu_\omega(Q_s)-C_su.
\]
The uniform slope gap $d_s$ therefore implies that each quotient bundle $\mathcal R_s$ is $\omega_u$-Higgs-stable for all sufficiently small $u>0$. Notice that $\mu_{\omega_u}(\mathcal R_s)\to\mu_\omega(Q_s)$, which is equal to  $\mu_{\omega}(\mathcal E_s / \mathcal E_{s-1})$. By setting
\[
 \lambda_s(u)=
 \frac{2\pi\,\mu_{\omega_u}(\mathcal R_s)}
      {\Vol(\widetilde Y,\omega_u)},
\]
we have $\lambda_s(u)=a_s+O(u)$ as $u \to 0$.

For each sufficiently small $u>0$, we have a smooth Higgs Hermitian--Einstein metric $k_{s,u}$ on each stable bundle $\mathcal R_s$ \cite{Simpson88}:
\[
 \ii\Lambda_{\omega_u}
 \bigl(\Theta_{k_{s,u}}
       +[\widetilde\theta_s,\widetilde\theta_s^{*k_{s,u}}]\bigr)
 =\lambda_s(u)\Id_{\mathcal R_s}.
\]

We use the standard block-scaling construction (cf.~\cite[Proposition 5.7]{LZ17}). Choose a smooth splitting $\widetilde\cB \simeq \bigoplus_s\mathcal R_s$. For $0<\tau<1$, define the direct-sum metric on $\widetilde\cB$:
\[
 H=H_{u,\tau}=\bigoplus_{s=1}^{m}\tau^{2(m-s)}k_{s,u}
\]
via the smooth isomorphism. Set
\[
 \Theta^{1,1}_{H,\widetilde\beta}
 =\Theta_H+[\widetilde\beta,\widetilde\beta^{*H}],\quad
 K_{u,H}=\ii\Lambda_{\omega_u}
          \Theta^{1,1}_{H,\widetilde\beta},\quad
 \Psi_{\lambda(u)}(H)=
 \sum_s\lambda_s(u)(P_s^H-P_{s-1}^H).
\]
For each fixed $u$, we can choose $\tau$ such that
\begin{equation}\label{eq:block-L1-input}
 \|K_{u,H}-\Psi_{\lambda(u)}(H)\|_{L^\infty}<u.
\end{equation}
In particular, $|K_{u,H}|_H$ is bounded independently of $u$ and the subsequently chosen $\tau$.

For fixed $u$, the Hitchin--Simpson Chern--Weil identity
\cite[proof of Proposition~3.4]{Simpson88} gives
\[
 \int_{\widetilde Y}
 \bigl(|\Theta^{1,1}_{H,\widetilde\beta}|_{\omega_u,H}^2+
       2|\partial_H\widetilde\beta|_{\omega_u,H}^2\bigr)dV_{\omega_u}
 =
 \int_{\widetilde Y}|K_{u,H}|_H^2dV_{\omega_u}
 +\mathcal T(u),
\]
where
\[
 \mathcal T(u)=4\pi^2\int_{\widetilde Y}
 \bigl(2c_2(\widetilde\cB)-c_1(\widetilde\cB)^2\bigr)
       \wedge\frac{\omega_u^{n-2}}{(n-2)!},  \  \mathcal T(u)=0 \  \mathrm{for} \  n=1.
\]
Since the volume and $\mathcal T(u)$ are bounded polynomials in $u$, we have
\begin{equation}\label{eq:block-L2-input}
 \|\Theta^{1,1}_{H,\widetilde\beta}\|_{L^2(\widetilde Y,\omega_u,H)}
 \leq C_0
\end{equation}
for one constant $C_0$ depending only on the fixed data.

Off the exceptional locus, locally choose a $\widetilde\eta$-unitary coframe which diagonalizes $\omega_0$:
\[
 \omega_0=\ii\sum_{j=1}^n\nu_j(x)\,
 e^j\wedge\overline{e^j},\quad 0\leq\nu_j(x)\leq M,
\]
where $M$ is the bound such that $\omega_0 \leq M \cdot \tilde \eta$. In this coframe,
\[
\omega_u=\ii\sum_{j=1}^n(\nu_j+u) \  e^j\wedge\overline{e^j}  \Rightarrow dV_{\omega_0}= \Big(\prod_{j=1}^n \nu_j\Big)dV_{\tilde \eta}, \ dV_{\omega_u}= \Big(\prod_{j=1}^n (\nu_j+u)\Big)dV_{\tilde \eta}.
\]
For any $\mathrm{End}(\widetilde \cB)$-valued $(1,1)$-form $A=\ii\sum_{j,k}A_{j\bar k}e^j\wedge\overline{e^k}$ off the exceptional locus, we have
\begin{equation*}
\begin{split}
\Lambda_{\omega_0}A-\Lambda_{\omega_u}A 
&=  \sum_j \Bigg(\frac{A_{j\bar j}}{\nu_j}-\frac{A_{j\bar j}}{\nu_j+u} \Bigg)
= \sum_j \frac{u}{\nu_j(\nu_j+u)}A_{j\bar j} \\
& \Rightarrow  |(\Lambda_{\omega_0}-\Lambda_{\omega_u})A|_H \cdot dV_{\omega_0} 
\leq \sum_j \frac{u}{\nu_j(\nu_j+u)}|A_{j\bar j}|  \Big(\prod_{k} \nu_k \Big) dV_{\tilde \eta} \\
& = u \sum_j \frac{\prod_{k \ne j} \nu_k}{\nu_j+u} |A_{j\bar j}| dV_{\tilde \eta}.
\end{split}
\end{equation*}
Notice that
\[
|A|_{\omega_u,H}^2= \sum_{j,k}\frac{|A_{j\bar k}|_H^2}{(\nu_j+u)(\nu_k+u)} \geq \sum_j \frac{|A_{j\bar j}|_H^2}{(\nu_j+u)^2},
\]
thus
\[
|(\Lambda_{\omega_0}-\Lambda_{\omega_u})A|_H  dV_{\omega_0} \leq u \sum_j \Big( \prod_{k \ne j} \nu_k \Big) |A|_{\omega_u,H} dV_{\tilde \eta} = c_u \cdot |A|_{\omega_u,H} dV_{\omega_u}
\]
for $c_u = u \frac{\sum_j \Big( \prod_{k \ne j} \nu_k \Big)}{\prod_j (\nu_j+u)}$. Applying the Cauchy--Schwarz inequality, we obtain:
\[
\|(\Lambda_{\omega_0}-\Lambda_{\omega_u})A\|_{L^1(\omega_0,H)} \leq \Big( \int_{\tilde Y} c_u^2 dV_{\omega_u} \Big)^\frac{1}{2}  \|A\|_{L^2(\omega_u,H)}.
\]

Finally, we can control the $L^2$-norm of $c_u$:
\begin{equation*}
      \begin{split}
c_u=\frac{u \sum_j \Big( \prod_{k \ne j} \nu_k \Big)}{\prod_{k} \nu_k + u \sum_j \Big( \prod_{k \ne j} \nu_k \Big) + \cdots} \leq 1; \\
\int_{\tilde Y} c_u dV_{\omega_u} = u \int_{\widetilde Y} \sum_j \Big( \prod_{k \ne j} \nu_k  \Big) dV_{\tilde \eta} = O(u)
      \end{split}
\end{equation*}
since $\nu_k \leq M$. Hence
\begin{equation}\label{eq:L1-contraction-comparison}
 \|(\Lambda_{\omega_0}-\Lambda_{\omega_u})A\|
       _{L^1(\omega_0,H)}
 \leq C_1\sqrt u\,\|A\|_{L^2(\omega_u,H)}.
\end{equation}
The contraction $\Lambda_{\omega_0}$ is defined off the exceptional locus of \(\pi\).

Write $\Psi_a(H)=\sum_s a_s(P_s^H-P_{s-1}^H)$ and
$K_{0,H}=\ii\Lambda_{\omega_0}\Theta^{1,1}_{H,\widetilde\beta}$.  Given $\delta>0$,  choose $u<\delta/(6\Vol(X,\omega))$  such that
\[
 C_1C_0\sqrt u<\frac{\delta}{6},\quad
 \Vol(X,\omega)\sum_s r_s|\lambda_s(u)-a_s|<\frac{\delta}{6}.
\]
(\ref{eq:block-L1-input}), (\ref{eq:block-L2-input}) and (\ref{eq:L1-contraction-comparison}) give
\[
 \begin{split}
 \|K_{0,H}-\Psi_a(H)\|_{L^1(\omega_0,H)}
 &\leq\|K_{0,H}-K_{u,H}\|_{L^1}
       +\|K_{u,H}-\Psi_{\lambda(u)}(H)\|_{L^1}
       +\|\Psi_{\lambda(u)}(H)-\Psi_a(H)\|_{L^1}\\
 &\leq C_1\sqrt u\,\|\Theta^{1,1}_{H,\widetilde\beta}\|_{L^2(\omega_u,H)}
       +\Vol(X,\omega)u
       +\Vol(X,\omega)\sum_s r_s|\lambda_s(u)-a_s|\\
 &<\frac{\delta}{2}.
 \end{split}
\]
The metric $H$ is now fixed.

Set $\widehat H:=\sigma^*\widehat h$.  If $\Sigma_f=\varnothing$, omit the cutoff
and set $H_r=H$.  Otherwise choose a cutoff function $\chi_r$ such that $\chi_r=1$ on the tubular neighborhood $N_{r/2}(\Sigma_f)$ and $\chi_r$ is supported in $N_r(\Sigma_f)$. We also require $\chi_r$ to satisfy
$|d\chi_r|_\omega\leq Cr^{-1}$ and
$|\partial\bar\partial\chi_r|_\omega\leq Cr^{-2}$.  Put
\[
 H_r=(1-\pi^*\chi_r)H+\pi^*\chi_r\widehat H,\quad
 V_r=\pi^{-1}(N_r(\Sigma_f)).
\]
$H_r$ is equal to $H$ outside $V_r$ and to $\widehat H$ near $\pi^{-1}(\Sigma_f)$, which includes all exceptional fibers of $\sigma$. This guarantees that $H_r$ can be pushed down to $\cB$. Moreover $\widehat H \sim H $ so $H_r \sim H $.

Choosing the same local $\widetilde{\eta}$-unitary coframe as before, we define $q_0 = dV_{\widetilde{\eta}}/dV_{\omega_0} = (\prod_{j=1}^n \nu_j)^{-1}$. For any $\mathrm{End}(\widetilde{\mathcal B})$-valued $(1,1)$-form $A$, $|\Lambda_{\omega_0}A|_{H_r} \leq \sqrt{n}\,|A|_{\omega_0,H_r}$, and $|A|_{\omega_0,H_r} \leq C_H q_0$ for fixed $A$. Furthermore, we have $|\Psi_a(H_r)|_{H_r} = \bigl(\sum_s r_s a_s^2\bigr)^{1/2} \leq C_H$. Now we estimate $\Theta^{1,1}_{H_r,\widetilde\beta}$:

Set $T=\widehat{H}^{-1} H_r$. Then
\[
\Theta_{H_r}-\Theta_{\widehat H}=\bar\partial(T^{-1}\partial_{\widehat H} T) = T^{-1}\bar\partial\partial_{\widehat H} T - T^{-1}\bar\partial T \cdot T^{-1}\partial_{\widehat H} T.
\]
Let $\varphi =\pi^* \chi_r$ and $S=\widehat{H}^{-1} H$. Then $T=(1-\varphi)S+\varphi \Id$. We have:

\begin{enumerate}
    \item $T^{-1}=H_r^{-1}\widehat H$ is bounded;
    
    \item 
    $\nabla_{\widehat H}T=d\varphi\cdot(\operatorname{Id}-S)
      +(1-\varphi)\nabla_{\widehat H}S =O(r^{-1})+O(q_0^{\frac12})$;
    \item
  
    $\bar\partial\partial_{\widehat H}T
    =\bar\partial\partial\varphi\cdot
      (\operatorname{Id}-S)
      +\partial_{\widehat H}\varphi\wedge\bar\partial S
    -\bar\partial\varphi\wedge
      \partial_{\widehat H}S
      +(1-\varphi)\bar\partial\partial_{\widehat H}S \\
    =O(r^{-2})+O(r^{-1} \cdot q_0^{\frac12})
    +O(q_0)
      =O(r^{-2}+q_0)$.
    
\end{enumerate}
Moreover, $\widetilde\beta$ and
$
\widetilde\beta^{*H_r}
=T^{-1}\widetilde\beta^{*\widehat H}T$
are both $O(q_0^{1/2})$. Hence $\Theta^{1,1}_{H_r,\widetilde\beta}=O(r^{-2}+q_0)$ and 
\begin{equation}\label{eq:L1-cutoff-residual}
 |K_{0,H_r}|_{H_r}+|\Psi_a(H_r)|_{H_r}
 \leq C_H(r^{-2}+1+q_0)\quad\hbox{almost everywhere on }V_r.
\end{equation}
The constant $C_H$ may depend on the already chosen $H,\widehat H$ but
is independent of $r$.

For the analytic subset $\Sigma_f$ of complex codimension at least 2, the volume of the tubular neighborhood $N_r(\Sigma_f)$ satisfies
\[
 \Vol_\omega(N_r(\Sigma_f))\leq C_{\rm tub}r^{4},
 \quad 0<r<r_0,
\]
where $C_{\rm tub}$ and $r_0$ depend only on $\Sigma_f$ and $\omega$. Since $q_0dV_{\omega_0}=dV_{\widetilde\eta}$, \eqref{eq:L1-cutoff-residual} yields
\[
 \begin{split}
 \int_{V_r}|K_{0,H_r}-\Psi_a(H_r)|_{H_r}dV_{\omega_0}
 &\leq C_H\Big(r^{-2}\Vol_\omega(N_r(\Sigma_f))
            +\Vol_\omega(N_r(\Sigma_f))
            +\Vol_{\widetilde\eta}(V_r)\Big)\\
 &\longrightarrow0.
 \end{split}
\]
Choose $r$ so that this integral is smaller
than $\delta/2$.  Outside $V_r$ the final metric is exactly $H$, with the
same norm.  Hence
\begin{equation}\label{before-pushdown}
 \|K_{0,H_r}-\Psi_a(H_r)\|_{L^1(\omega_0,H_r)}<\delta.
\end{equation}
Since $\sigma$ is biholomorphic outside $\pi^{-1}(\Sigma_f)$ and $H_r = \sigma^*\widehat h$ on $V_{r/2}$, we can smoothly push down $H_r$ to obtain a Hermitian metric $\widehat h_\delta$ on $\mathcal B$, agreeing with \(\widehat h\) on \(\rho^{-1}(N_{r/3}(\Sigma_f))\), such that $\sigma^*\widehat h_\delta=H_r$. By omitting analytic subsets, we can also push down the endomorphisms $K_{0,H_r}$ and $\Psi_a(H_r)$ to get $K_\delta$ and $\sum_s a_s(P_s-P_{s-1})$, and hence \eqref{before-pushdown}. Therefore \eqref{eq:approximate-critical} holds.

For \eqref{eq:regularized-residual-convergence}, keep $\widehat h_\delta$
fixed. On $\rho^{-1}(X_f)$, put
\[
 \Psi_\delta
 =\sum_s a_s(P_{s,\delta}-P_{s-1,\delta}) = \rho^*\Big(\sum_s a_s(P_s-P_{s-1})\Big).
\]
Clearly the densities converge pointwise on $\rho^{-1}(X_f)$:
\[
 |K_{\delta,\epsilon}-\Psi_\delta|_{\widehat h_\delta}
 dV_{\omega_\epsilon}
 \longrightarrow
 \rho^*\big(
|K_\delta-\sum_s a_s(P_s-P_{s-1})|_{h_\delta}
 dV_\omega \big).
\]
It remains to find an integrable dominating density.

Locally choose a $\eta$-unitary coframe which diagonalizes $\rho^*\omega$:
\[
\rho^*\omega=\ii\sum_{j=1}^n L_j^{\rho}e^j\wedge\overline{e^j} \Rightarrow dV_{\omega_\epsilon}
               =\prod_j(L_j^\rho+\epsilon) dV_\eta
\]
with $0\leq L_j^\rho\leq M_\rho$. The same argument as in the proof of \eqref{eq:L1-contraction-comparison} gives

\[
|\Lambda_{\omega_\epsilon}
\Theta^{1,1}_{\widehat h_\delta,\beta}|_{\widehat h_\delta}dV_{\omega_\epsilon}
 \leq C_\delta\big(\sum_j\prod_{k\ne j}(L_k^\rho+\epsilon)\big)dV_\eta
 \leq C_\delta n(M_\rho+1)^{n-1} dV_\eta.
\]
Moreover
$|\Psi_\delta|_{\widehat h_\delta}
 =(\sum_s r_sa_s^2)^{1/2}$ on $\rho^{-1}(X_f)$.  Thus, for fixed $\delta$,
\[
 |K_{\delta,\epsilon}-\Psi_\delta|_{\widehat h_\delta}dV_{\omega_\epsilon}
 \leq C'_\delta dV_\eta
 \quad\text{almost everywhere on }Y,\quad 0<\epsilon\leq1.
\]
\eqref{eq:regularized-residual-convergence} is then a consequence of the dominated convergence theorem.
\end{proof}

\begin{proposition}\label{prop:HN-vector}
For every integer $q$ with $1\leq q\leq R$, $\mu_q=b_q$.

\end{proposition}

\begin{proof} 
Instead of establishing the continuous dependence of the heat flow on initial metrics over the non-compact $X_0$ (cf. \cite[Lemma 5.2]{LZZ18}), we use the Higgs Donaldson functional to compare heat flows on the regularized Higgs bundle \(\mathcal B\) directly. Write
\[
 V_\epsilon=\Vol(Y,\omega_\epsilon),\quad
 K_{\epsilon,H}=\ii\Lambda_{\omega_\epsilon}
       \bigl(\Theta_H+[\beta,\beta^{*H}]\bigr),\quad
 \Phi_{\epsilon,H}=K_{\epsilon,H}-\lambda_\epsilon\Id,
\]
where the constant $\lambda_\epsilon=2\pi\mu_{\omega_\epsilon}(\cB)/V_\epsilon$. Section~\ref{sec:heat-kernel-estimates}, in particular 
\eqref{eq:initial-phi-bound} and \eqref{eq:uniform-heat-kernel}, provides
\[
 \begin{gathered}
 V_\epsilon\leq V_+,\quad
 V_\epsilon\longrightarrow\Vol(X,\omega),\quad
 \lambda_\epsilon\longrightarrow\lambda,\\
 M_1:=\sup_{0<\epsilon\leq1}
 \int_Y|\Phi_{\epsilon,\widehat h}|_{\widehat h}
 dV_{\omega_\epsilon}<\infty,\quad
 C_{\rm heat}:=\sup_{0<\epsilon\leq1}\sup_{x,y\in Y}
 H_\epsilon(x,y,1)<\infty.
 \end{gathered}
\]
Notice that the heat kernel $H_\epsilon$ is defined with respect to the volume form $dV_{\omega_\epsilon}$.

Fix $0<\delta<1$, and let $\widehat h_\delta$ be the metric supplied by
Lemma~\ref{lem:approximate-critical}. Retain the notation $\Psi_\delta=\sum_s a_s(P_{s,\delta}-P_{s-1,\delta})$.
By \eqref{eq:approximate-critical} and
\eqref{eq:regularized-residual-convergence}, there is
$\epsilon_\delta>0$ such that
\[
 Q_{\delta,\epsilon}
 =K_{\delta,\epsilon}
  -\Psi_\delta,\quad
 \|Q_{\delta,\epsilon}\|_{L^1}<\delta
 \quad(0<\epsilon\leq\epsilon_\delta).
\]
Recall that the eigenvalues of
$\Psi_\delta$ are exactly
$\mu_1,\ldots,\mu_R$ almost everywhere.  Define the constants
\[
 A_{\rm HN}=\sup_{0<\epsilon\leq1}
      \left(\sum_q(\mu_q-\lambda_\epsilon)^2\right)^\frac{1}{2}=\sup_{0<\epsilon\leq1}\Big| \Psi_\delta - \lambda_{\epsilon} \Id   \Big|_{\widehat h_\delta},
 \quad M_2=1+V_+A_{\rm HN}.
\]
They are finite and independent of $\delta$.  It follows that 
\begin{equation}\label{eq:initial-phi-bound-delta}
\Phi_{\epsilon,\widehat h_\delta}
 =K_{\delta,\epsilon}-\lambda_\epsilon\Id
 =Q_{\delta,\epsilon}+(\Psi_\delta-\lambda_\epsilon\Id)   \Rightarrow
\|\Phi_{\epsilon,\widehat h_\delta}\|_{L^1}\leq M_2
\end{equation}
for $0<\epsilon\leq\epsilon_\delta$.

On the smooth Higgs bundle $\mathcal B$, consider the following two heat flows:
\[
\left\{
\begin{aligned}
h_\epsilon^{-1}\partial_t h_\epsilon &= -2\Phi_{\epsilon,h_\epsilon}, & h_\epsilon(0)&=\widehat h,\\
k_{\delta,\epsilon}^{-1}\partial_t k_{\delta,\epsilon} &= -2\Phi_{\epsilon,k_{\delta,\epsilon}}, & k_{\delta,\epsilon}(0)&=\widehat h_\delta.
\end{aligned}
\right.
\]

For fixed $\epsilon$, the long-time existence for both flows is guaranteed by \cite[Proposition~6.6]{Simpson88}. Define their energies by
\[
 E_\epsilon^h(t)=
   \int_Y|\Phi_{\epsilon,h_\epsilon(t)}|^2_{h_\epsilon(t)}
         dV_{\omega_\epsilon},\quad
 E_{\delta,\epsilon}^k(t)=
   \int_Y|\Phi_{\epsilon,k_{\delta,\epsilon}(t)}|^2_{k_{\delta,\epsilon}(t)}
         dV_{\omega_\epsilon}.
\]
Direct calculation gives $(\partial_t-\Delta_{\omega_\epsilon})|\Phi|^2=-2|\nabla^{\mathrm{HS}}\Phi|^2$, and Kato's inequality gives $(\partial_t-\Delta_{\omega_\epsilon})|\Phi|\leq0$ in the weak sense.
Integrating over $Y$ and using the kernel bound $C_{\rm heat}$, we obtain
\begin{equation}\label{eq:energy-bound}
 \|\Phi_{\epsilon,h_\epsilon(t)}\|_{L^1}\leq M_1,\quad
 \|\Phi_{\epsilon,k_{\delta,\epsilon}(t)}\|_{L^1}
       \leq M_2,\quad
 \|\Phi_{\epsilon,k_{\delta,\epsilon}(1)}\|_{L^\infty}
       \leq C_{\rm heat}M_2
\end{equation}
for all $t\geq0$. Fix $L>C_{\rm heat}M_2$, independent of $\delta$ and $\epsilon$. Define
\[
 f_L(x)=
 \begin{cases}
 x^2,&|x|\leq L,\\
 2L|x|-L^2,&|x|>L.
 \end{cases}
\]
This convex $C^1$ function satisfies $|f_L'|\leq2L$.  For two self-adjoint
endomorphisms $A,B$ in the same Hermitian metric, the real-valued function
$\tr f_L(B+s(A-B))$ is also $C^1$ in $s$, and its derivative is $\tr(f_L'(B+s(A-B))(A-B))$. 
 \[
 |\tr f_L(A)-\tr f_L(B)|\leq2L\sqrt R\,|A-B|.
\]
Here $R=\operatorname{rank}(\mathcal B)$ and the norm on the right is the Hilbert--Schmidt norm.

For a smooth convex function $f$, let $f'^{[1]}(x,y)$ be $(f'(x)-f'(y))/(x-y)$, with value $f''(x)$ at $x=y$. The heat flow of the metric $H$ gives the flow of  $\Phi$:
\[
 \partial_t\Phi=-2\ii\Lambda_{\omega_\epsilon}
 \bigl(\bar\partial \partial_H\Phi+
           [\beta,[\beta^{*H},\Phi]]\bigr).
\]
In a local $H$-unitary frame with
$\Phi=\operatorname{diag}(\nu_1,\ldots,\nu_R)$, the calculation in \cite[Lemma 6.1]{Simpson88} gives the monotonicity of the spectral functional:
\[
 \frac{d}{dt}\int_Y\tr f(\Phi)dV_{\omega_\epsilon}
 \overset{\mathrm{loc}}{=}
-2\int_Y\sum_{a,b}f'^{[1]}(\nu_a,\nu_b)
 \bigl(|(\bar\partial \Phi)_{ab}|^2
+|[\beta,\Phi]_{ab}|^2\bigr)dV_{\omega_\epsilon}\leq0.
\]

By approximating $f_L$ uniformly with smooth convex functions, we obtain the monotonicity of its spectral functional. At $t=1$, the choice of $L$ makes $f_L(x)=x^2$ for $x=\Phi_{\epsilon,k_{\delta,\epsilon}(1)}$; at $t=0$, recalling \eqref{eq:initial-phi-bound-delta} with $\|Q_{\delta,\epsilon}\|_{L^1}<\delta$, we have
\[
\int_Y \left|
\operatorname{tr}f_L(\Phi_{\epsilon,\widehat h_\delta})
-
\operatorname{tr}f_L(\Psi_\delta-\lambda_\epsilon \Id)
\right| dV_{\omega_\epsilon}
\leq
2L\sqrt R\,\|Q_{\delta,\epsilon}\|_{L^1} \leq
2L\sqrt R\,\delta ,
\]
together with 
\[
\int_Y\operatorname{tr}
f_L(\Psi_\delta-\lambda_\epsilon \Id)\,dV_{\omega_\epsilon}
=
V_\epsilon\sum_qf_L(\mu_q-\lambda_\epsilon)
\leq
V_\epsilon\sum_q(\mu_q-\lambda_\epsilon)^2.
\]
Hence the energy at $t=1$ is bounded by

\begin{equation}\label{eq:initial-spectral-bound}
E_{\delta,\epsilon}^k(1) \leq \int_Y\tr f_L(\Phi_{\epsilon,\widehat h_\delta})dV_{\omega_\epsilon}
\leq V_\epsilon\sum_q(\mu_q-\lambda_\epsilon)^2
       +2L\sqrt R\,\delta.
\end{equation}
Taking \(f(x)=x^2\) also shows that both energies are nonincreasing.

On the Higgs bundle $(Y,\omega_\epsilon,\cB,\beta)$, the Higgs
Donaldson functional is defined by
\[
 \mathcal M_\epsilon(H_0,H_1)
  =\int_0^1\int_Y\tr(H_t^{-1}\partial_tH_t\,\Phi_{\epsilon,H_t})
                 dV_{\omega_\epsilon}dt.
\]

\cite[Proposition~5.1 and its proof]{Simpson88} already yields the independence of the choice of path $H_t$, the convexity along the geodesic $H_t=Ke^{tS}$ with $S=S^{*K}$, and the cocycle identity
\[
\mathcal M_\epsilon(H_0,H_1)
+\mathcal M_\epsilon(H_1,H_2)
=\mathcal M_\epsilon(H_0,H_2).
\]
Note that our definition presented here is closer to the original one in \cite{Donaldson87} (without Higgs fields) than the one in \cite{Simpson88}, but they are equivalent owing to path-independence.

Consider the relative endomorphism $q_{\delta,\epsilon}=h_\epsilon^{-1}k_{\delta,\epsilon}$. The same calculation as in \cite[Proposition~13]{Donaldson85}, using \cite[Lemma~3.1]{Simpson88}, shows that

\[
 (\Delta_{\omega_\epsilon}-\partial_t)
 \tr(q_{\delta,\epsilon}+q_{\delta,\epsilon}^{-1}-2\Id)\geq0.
\]
By the maximum principle, $\tr(q_{\delta,\epsilon}+q_{\delta,\epsilon}^{-1}-2\Id)$ is  bounded uniformly for $t\geq0$, which implies
\begin{equation}\label{eq:comparison-metric-equivalence}
 C_\delta^{-1}h_\epsilon(t)
 \leq k_{\delta,\epsilon}(t)\leq C_\delta h_\epsilon(t),
 \quad t\geq0,
\end{equation}
where $C_\delta\geq1$ is independent of $\epsilon$ and $t$.
At a fixed time $t$ put
$S=\log(k_{\delta,\epsilon}^{-1}h_\epsilon)$ and
$H_s=k_{\delta,\epsilon}e^{sS}$.
The endomorphism $S$ is self-adjoint along the  geodesic and satisfies
$|S|_{H_s}=\tr(S^2)^{1/2} \leq \sqrt R\log C_\delta$.
For the convex function $g$  defined by
\[
g(s)=\mathcal M_\epsilon(k_{\delta,\epsilon},H_s) \Rightarrow g'(s)=\int_Y\tr(S\Phi_{\epsilon,H_s})dV_{\omega_\epsilon},
\]
the inequality $g'(0) \leq g(1)-g(0) \leq g'(1)$ implies
\[
 \int_Y\tr(S\Phi_{\epsilon,k_{\delta,\epsilon}})dV_{\omega_\epsilon}
 \leq\mathcal M_\epsilon(k_{\delta,\epsilon},h_\epsilon)
 \leq\int_Y\tr(S\Phi_{\epsilon,h_\epsilon})dV_{\omega_\epsilon}.
 \]
Now \eqref{eq:energy-bound} gives
\[
 \left|\mathcal M_\epsilon
   (k_{\delta,\epsilon}(t),h_\epsilon(t))\right|
 \leq A_\delta,\quad
 A_\delta=\sqrt R\log C_\delta\max\{M_1,M_2\}.
\]
The bound for the Donaldson functional is uniform for
$t\geq0$ and $0<\epsilon\leq\epsilon_\delta$.

Choose any fixed reference metric $K$. The cocycle identity gives
\[
\mathcal M_\epsilon(k_{\delta,\epsilon}(t),h_\epsilon(t))
=\mathcal M_\epsilon(K,h_\epsilon(t)) -\mathcal M_\epsilon(K,k_{\delta,\epsilon}(t)).
\]
Differentiating the above identity and using the heat flow equation $H^{-1}\partial_tH=-2\Phi$ gives
\begin{equation*}
 \frac{d}{dt}\mathcal M_\epsilon
   (k_{\delta,\epsilon}(t),h_\epsilon(t))
 =-2E_\epsilon^h(t)+2E_{\delta,\epsilon}^k(t).
\end{equation*}
Integrating from $1$ to $T>1$ and using the energy monotonicity together with the bound $A_\delta$ of $\mathcal M_\epsilon(k_{\delta,\epsilon}(t),h_\epsilon(t))$, we get
\[
 (T-1)E_\epsilon^h(T)
\leq\int_1^T E_\epsilon^h(t)\,dt
\leq\int_1^T E_{\delta,\epsilon}^k(t)\,dt+A_\delta
\leq(T-1)E_{\delta,\epsilon}^k(1)+A_\delta.
\]
Combining this with \eqref{eq:initial-spectral-bound}, we prove that
\begin{equation}\label{eq:compact-comparison-energy}
 E_\epsilon^h(T)
 \leq V_\epsilon\sum_q(\mu_q-\lambda_\epsilon)^2
        +2L\sqrt R\,\delta+\frac{A_\delta}{T-1}.
\end{equation}

Keep $\delta,T$ fixed and let $\epsilon=\epsilon_j\downarrow0$ as in Proposition~\ref{prop:positive-flow}. The local smooth convergence of $h_{\epsilon_j}(t)$ shows that, by Fatou's lemma,
\[
 E^h(T):=\int_{X_0}|\Phi(T)|_{h(T)}^2\ddV
\leq\liminf_{j\to\infty}E_{\epsilon_j}^h(T).
\]
Then \eqref{eq:compact-comparison-energy} yields
\[
 E^h(T)\leq E_{\rm HN}+2L\sqrt R\,\delta+\frac{A_\delta}{T-1},
 \quad E_{\rm HN}=\Vol(X,\omega)\sum_q(\mu_q-\lambda)^2.
\]
First let $T\to\infty$ for fixed $\delta$, and then let $\delta\downarrow0$.
Since $L$ is independent of $\delta$, we get $\limsup_{T\to\infty}E^h(T)\leq E_{\rm HN}$. 

For the times $t_i \to \infty$ chosen in Proposition~\ref{prop:slice-compactness}, we have the following convergence in the local unitary gauges $\{a_{i,\alpha}\}$:
\[
{a_{i,\alpha}} \Psi(t_i) {a_{i,\alpha}}^{-1} \rightarrow \Psi_\infty \Rightarrow |\Phi(t_i)|_{h(t_i)}^2 \rightarrow |\Psi_\infty|_{h_\infty}^2 
\]
locally uniformly on $U=X \setminus Z$, with $|\Psi_\infty|_{h_\infty}^2 = \tr (\Psi_\infty^2)=\sum_q(b_q-\lambda)^2$. Therefore
\begin{equation}\label{eq:reverse-power}
 \Vol(X,\omega)\sum_q(b_q-\lambda)^2
 \leq\liminf_{i\to\infty}E^h(t_i)
 \leq\Vol(X,\omega)\sum_q(\mu_q-\lambda)^2
\end{equation}
by Fatou's lemma again.
Finally, Lemma~\ref{lem:polygon} provides
\[
 \sum_{q=1}^R(b_q-\mu_q)=0,\quad
 \sum_{q=1}^k(b_q-\mu_q)\geq0\quad(1\leq k<R).
\]
Direct calculation shows that
\[
 \begin{split}
 \sum_q(b_q-\lambda)^2-\sum_q(\mu_q-\lambda)^2
 &=\sum_q(b_q-\mu_q)^2
   +2\sum_q(\mu_q-\lambda)(b_q-\mu_q)\\
 &=\sum_q(b_q-\mu_q)^2
   +2\sum_{k=1}^{R-1}(\mu_k-\mu_{k+1})
       \sum_{q=1}^k(b_q-\mu_q)\\
 &\geq\sum_q(b_q-\mu_q)^2 \geq 0.
 \end{split}
\]
However, by \eqref{eq:reverse-power}, the left side is nonpositive, so $b_q=\mu_q$
for every $q$.
\end{proof}

\section{Proof of the main theorem}
\label{sec:residual-exhaustion}

We prove the theorem by successively splitting off the stable factors $Q_s=(\mathcal E_s/\mathcal E_{s-1})^{\vee\vee}$ from $(\mathcal E_{\infty},\theta_{\infty})$. Our argument proceeds by induction, with the following inductive hypotheses:

\textbf{Inductive hypotheses} \begin{itshape} Recall the bad set $Z$ in Proposition \ref{prop:slice-compactness} and $\Sigma_f$ in Lemma \ref{lem:approximate-critical}. Put $U_f=X\setminus(Z\cup\Sigma_f)$. When necessary, we pass to a subsequence of $(A(t_i),\theta(t_i))$ in Proposition \ref{prop:slice-compactness} and denote its times by $t_j$.

 \begin{enumerate}
    \item At stage $s$, the first $s-1$ summands have been split:
\[
 (\cE_{\infty},\theta_{\infty})
 \simeq \left(\bigoplus_{q=1}^{s-1}Q_q\right)\oplus T_{s-1},
\]
where $T_{s-1}$ is a direct sum of reflexive polystable Higgs sheaves.

    \item Let $g_j = g(t_j)$. Define $g_j(\mathcal E_s)$ to be the bundle obtained by the action of the complex gauge transformation $g_j$ on $\mathcal E_s$ over $U_f$, and $P_{s-1,j}: g_j(\cE) \to g_j(\cE_{s-1})$ to be the  $h_0$-orthogonal projection. Then $P_{s-1,j}\longrightarrow P_{s-1,\infty}$ in the sense of $C^{\infty}_{\mathrm{loc}}$ in local unitary gauges $\{a_{j,\alpha}\}$ $($Prop. \ref{prop:slice-compactness}$)$.

    \item Define the bundle $T_{s-1,j}=g_j(\mathcal E)/g_j(\mathcal E_{s-1})$ over $U_f$. Then $T_{s-1,j}\to T_{s-1}$, also in $C^\infty_{\mathrm{loc}}$, with respect to the local unitary gauges $\{a_{j,\alpha}\}$. Here the convergence includes holomorphic structures and Higgs fields.

\end{enumerate}
\end{itshape}
The exact forms of $P_{s,\infty}$ and $T_{s}$ are shown in the proof of Lemma~\ref{lem:quotient-propagation}.

At stage $s$, recall $\pi=\rho\circ\sigma$ in Lemma \ref{lem:approximate-critical} and let
 $S=\widetilde{\mathcal E}_s/\widetilde{\mathcal E}_{s-1}$.
The composition
\[
\mathcal E_s / \mathcal E_{s-1} \to \mathcal E / \mathcal E_{s-1} \simeq g_j(\mathcal E)/g_j(\mathcal E_{s-1})
\]
over $U_f$ induces an injective Higgs morphism
\begin{equation}\label{eq:physical-index-map}
 u_{s,j}:Q_s|_{U_f}\longrightarrow T_{s-1,j}|_{U_f}.
\end{equation}

Fix a smooth Hermitian metric $h_S$ on $S$.
Via the isomorphism $S|_{\pi^{-1}(U_f)}\simeq\pi^*(Q_s|_{U_f})$ induced by the biholomorphic map $\pi:\pi^{-1}(U_f)\to U_f$,
$h_S$ descends to a metric $h_Q$ on $Q_s|_{U_f}$.
 $T_{s-1,j}$ is endowed with the quotient metric induced by $h_0$, which is the same as the quotient metric induced by $h(t_j)$ on $\cE/\cE_{s-1}$ due to
 \eqref{eq:horizontal-isometry}.
Write $|\cdot|_{0,j}$ for the $\mathrm{Hom}(Q_s, T_{s-1,j})$-norm defined by $h_Q$ and the quotient metric on $T_{s-1,j}$.

\begin{proposition}
\label{prop:normalized-maps}
Under the inductive hypotheses at stage $s$, there exist positive constants $c_{s,j}$ and a further subsequence of $\{t_j\}$ such that the maps
$f_{s,j}=c_{s,j}u_{s,j}$ converge locally smoothly on $U_f$.
The limiting Higgs morphism is nonzero and extends uniquely to
\begin{equation}\label{eq:limiting-residual-map}
 f_{s,\infty}:Q_s\longrightarrow T_{s-1}.
\end{equation}
\end{proposition}

\begin{proof}
Put $T=\widetilde\cB/\widetilde\cE_{s-1}$, and endow $T$ with the quotient metric induced by $\sigma^*h_{\epsilon_k}(t_j)$.
The natural bundle embedding gives a fiberwise injective bundle map
$\iota_s:S\to T$. (Note that fiberwise injectivity of a bundle map is a stronger condition than injectivity as a sheaf morphism.) Firstly, as in Section~\ref{sec:slice-compactness}, we consider the complex gauge transformation $g_{\epsilon_k}(t)$ such that
\[
\partial_t g_{\epsilon_k}=-g_{\epsilon_k}\Phi_{\epsilon_k},\quad g_{\epsilon_k}(0)=\Id,
\]
which implies that
 \[
 h_{\epsilon_k}(t)(u,v)=h_{\epsilon_k}(0)(g_{\epsilon_k}(t)u,g_{\epsilon_k}(t)v).
 \]
  By pulling back to $\widetilde Y$, we set $G_{j,k}=\sigma^*g_{\epsilon_k}(t_j)$ and $T_{s-1,j,k}=G_{j,k}(\widetilde{\mathcal B}) / G_{j,k}(\widetilde{\mathcal E}_{s-1})$, endowed with the quotient metric induced by $\sigma^*h_{\epsilon_k}(0)$. We can get an injective Higgs morphism $u_{s,j,k}$, constructed in the same way as $u_{s,j}$ from \eqref{eq:physical-index-map}. Its squared $\mathrm{Hom}(S,T_{s-1,j,k})$-norm $|u_{s,j,k}|_0^2$ equals the squared $\mathrm{Hom}(S,T)$-norm $|\iota_s|^2$, which is a smooth, strictly positive function on $\widetilde Y$.

Fix a Kähler form $\widetilde\eta$ with
$\Vol(\widetilde Y,\widetilde\eta)=1$, and consider the two-parameter family:
\[
 \omega_{k,a}=\pi^*\omega+\epsilon_k\sigma^*\eta+a\widetilde\eta,
 \quad 0<\epsilon_k,a\leq1.
\]
Let $H_{k,a}(x,y,t)$ be the heat kernel for $\Delta_{\omega_{k,a}}=2\ii\Lambda_{\omega_{k,a}}\partial\bar\partial$. An argument similar to that in the proof of \cite[Proposition~2]{BS94} yields the uniform bound at $t=1$:
\begin{equation}\label{eq:scalar-aux-kernel} 
 \sup_{k,a,x,y}H_{k,a}(x,y,1)\leq H_3.
\end{equation}
Here $H_3$ is independent of $k,a$. Note that 
$\omega_{k,a}\leq \nu:=\pi^*\omega+\sigma^*\eta+\widetilde\eta \Rightarrow dV_{\omega_{k,a}}\leq dV_{\nu}$. 

We define
\[
 K_B=\ii\Lambda_{\omega_{k,a}}
   \bigl(\Theta_{h_B,\theta_B}^{1,1 }\bigr),\quad
 b_B=(1+|K_B|_{h_B}^2)^{1/2},\quad b=b_S+b_T.
\]
The goal is to prove the uniform boundedness of $\int_{\widetilde Y}b\,dV_{\omega_{k,a}}$. The case $b_S$ is easy since the metric $h_S$ is fixed. The wedge identity provides:
\[ K_S dV_{\omega_{k,a}}
= \ii \bigl(\Theta_{h_S,\theta_S}^{1,1 }\bigr) \wedge \frac{\omega_{k,a}^{n-1}}{(n-1)!} \Rightarrow \int_{\widetilde Y} |K_S|_{h_S}   \,dV_{\omega_{k,a}} \leq C_S \cdot \Vol(\widetilde Y,\nu).
\]

For $b_T$ induced by the $t_j$-dependent metric $\sigma^*h_{\epsilon_k}(t_j)$, recall that
\[
\begin{gathered}
\int_Y|\Phi_{\epsilon_k}(t)|\,dV_{\omega_{\epsilon_k}}
\leq\int_Y|\Phi_{\epsilon_k}(0)|\,dV_{\omega_{\epsilon_k}}, \quad \lambda_{\epsilon_k} \to \lambda \\
\Rightarrow \int_Y|K_{\mathcal B}(t)|\,dV_{\omega_{\epsilon_k}}
\leq M_1 + (\mathrm{rank}\cB)|\lambda_{\epsilon_k}|\Vol(Y,\omega_{\epsilon_k}) \leq M_\cB,
\end{gathered}
\]
where $M_1$ and $K_{\mathcal B}=K_{\epsilon_k, h_{\epsilon_k}}$ are from Proposition~\ref{prop:HN-vector}. Now let
\[
G_a \cdot dV_{\nu} := \ii\Lambda_{\omega_{k,a}}
   \bigl(\sigma^* \Theta_{h_{\epsilon_k} (t_j),\beta}^{1,1 }\bigr) dV_{\omega_{k,a}} =  \ii \sigma^* \Theta_{h_{\epsilon_k} (t_j),\beta}^{1,1 }\wedge \frac{\omega_{k,a}^{n-1}}{(n-1)!}.
\]
By expanding $\omega_{k,a}^{n-1}$, we have the smooth convergence
\[
\ii \sigma^* \Theta_{h_{\epsilon_k} (t_j),\beta}^{1,1 }\wedge \frac{\omega_{k,a}^{n-1}}{(n-1)!} \longrightarrow
\ii \sigma^* \Big(\Theta_{h_{\epsilon_k} (t_j),\beta}^{1,1 }\wedge \frac{\omega_{\epsilon_k}^{n-1}}{(n-1)!} \Big) = \sigma^* \big(K_{\epsilon_k}(t_j) dV_{\omega_{\epsilon_k}}\big) 
\]
as $a\downarrow0$. Taking absolute values and integrating over $\widetilde Y$, we get
\[
 \int_{\widetilde Y}
 |\ii\Lambda_{\omega_{k,a}}\sigma^*
   \Theta^{1,1}_{h_{\epsilon_k}(t_j),\beta}|\,dV_{\omega_{k,a}}
 \longrightarrow
 \int_Y|K_{\mathcal B}(t_j)|\,dV_{\omega_{\epsilon_k}}  \leq M_{\cB}.
\]
Thus $||K_{\widetilde \cB}||_{L^1} \leq M_{\cB}+1$ for sufficiently small $a$. Now we apply the Gauss--Codazzi formula (cf. \cite[Chapter~V, Theorem 14.5]{Demailly12}) with Higgs fields to the bundle exact sequence
\[
0\xrightarrow{} \widetilde\cE_{s-1}
\xrightarrow{j} \widetilde\cB
\xrightarrow{g} T
\xrightarrow{} 0
\]
with metric $\sigma^*h_{\epsilon_k}(t_j)$ on $\widetilde \cB$, which gives
\[
 K_T=gK_{\widetilde\cB}g^*+R_T,\quad R_T\geq0.
\]
Using the facts that \(\|g\|_{\mathrm{op}}\le1\), $|R_T|\leq\operatorname{tr} R_T$ and $\int_{\widetilde Y} \operatorname{tr} K_T \,dV_{\omega_{k,a}}=2\pi\deg_{\omega_{k,a}} T$, we obtain:
\[
 \begin{aligned}
 ||K_T||_{L^1} & \leq ||g K_{\widetilde\cB}g^*||_{L^1} + \int_{\widetilde Y} \tr R_T dV_{\omega_{k,a}} \\ & \leq ||K_{\widetilde\cB}||_{L^1} + \int_{\widetilde Y} \operatorname{tr} K_T \,dV_{\omega_{k,a}} - \int_{\widetilde Y} \tr (g K_{\widetilde\cB}g^*) \,dV_{\omega_{k,a}} \\ &
 \leq ||K_{\widetilde\cB}||_{L^1} + 2\pi \deg_{\omega_{k,a}} T + \int_{\widetilde Y} |\tr (g K_{\widetilde\cB}g^*)| \,dV_{\omega_{k,a}}  \\
 & \leq (1+ \sqrt{\mathrm{rank} (T)}) ||K_{\widetilde\cB}||_{L^1} + 2\pi \deg_{\omega_{k,a}} T.
 \end{aligned}
\]
Since $b_T\leq1+|K_T|$, $b_T$ is uniformly bounded in $L^1$. Let $M_3$ be a uniform bound for $||b||_{L^1}$. Let
$A=H_3e^{2H_3M_3}$.

Write $q=|\iota_s|^2$. Applying the Bochner identity for the section $\iota_s$ of $W=\operatorname{Hom} (S,T)$ and using the Higgs fields, we have:
\[
\begin{aligned}
\Delta_{\omega_{k,a}}q &   = 
2|\partial_W \iota_s|^2
-
2\operatorname{Re}
\langle
\sqrt{-1}\Lambda_{\omega_{k,a}}\Theta_W(\iota_s),
\iota_s
\rangle \\
 & = 2|\partial_W \iota_s|^2
+2|\theta_T^{*h_T}\iota_s
-\iota_s\theta_S^{*h_S}|^2
-2\operatorname{Re}
\langle K_T\iota_s-\iota_sK_S,\iota_s\rangle \\
      & \geq 2|\partial_W \iota_s|^2 - 2bq.
\end{aligned}
\]
Moreover, the Cauchy--Schwarz inequality gives (note that $\bar \partial \iota_s=0$)
\[
|dq|^2 =
|\langle\nabla_W\iota_s,\iota_s\rangle
+
\langle\iota_s,\nabla_W\iota_s\rangle|^2 \leq 2|\langle \partial_W\iota_s,\iota_s\rangle|^2 \leq 2q|\partial_W\iota_s|^2.
\]
Combining the two inequalities above, we obtain the logarithmic estimate
\begin{equation}\label{eq:log-Hom-bochner}
 \Delta_{\omega_{k,a}}\log q\geq-2b.
\end{equation}

Write $\Delta=\Delta_{\omega_{k,a}}$ and set
$\psi=\int_0^1e^{t\Delta}b\,dt$. The basic properties of the heat equation give
\[
 \psi\geq0,\quad
 \Delta\psi=e^{\Delta}b-b= \int_{\widetilde Y} H_{k,a}(\cdot,y,1) b(y) \,dV_{\omega_{k,a}}(y)-b,\quad
 \int_{\widetilde Y}\psi\,dV_{\omega_{k,a}}
 =\int_{\widetilde Y}b\,dV_{\omega_{k,a}}\leq M_3.
\]
Let $w=e^{-2\psi}q$. Then
$\Delta\log w\geq-2e^{\Delta}b\geq-2H_3M_3 \Rightarrow \Delta w\geq-cw$, where $c=2H_3M_3$.
Thus $\partial_t\big(e^{ct}e^{t\Delta}w\big)\geq0$. Integrating from $t=0$ to $t=1$ yields

\begin{equation}\label{eq:scalar-weight-density}
 \sup_{\widetilde Y}w \leq ||e^ce^{\Delta}w||_{L^\infty(\widetilde Y)}   \leq A\int_{\widetilde Y}w\,dV_{\omega_{k,a}}.
\end{equation}

Fix $L\subset O\Subset\pi^{-1}(U_f)$, where $L$ is compact and $O$ is open. Since $\pi^* \omega>0$ on $\pi^{-1}(U_f)$, we have $\mathrm{sup}_O |K_S| \leq C_O^S$. For $K_T$, recall the local smooth convergence from Proposition \ref{prop:slice-compactness} and the inductive hypothesis $P_{s-1,j} \to P_{s-1,\infty}$. By the preceding Gauss--Codazzi formula, we obtain uniform control on the quotient Chern curvatures and quotient Higgs fields over $O$. In particular
\[
\sup_O |K_T| \leq C_O^T \Rightarrow b \leq C_O := 2+C_O^S+C_O^T
\]
Now that $\Delta \psi \geq -b \geq -C_O$, the standard Moser iteration with a ball covering argument gives
\begin{equation}\label{eq:local-scalar-weight}
 \sup_L\psi
 \leq C_{L,O}\left(\int_O\psi\,dV_{\omega_{k,a}}+C_O\right)
 \leq C_{L,O}(M_3+C_O)=:B_L.
\end{equation}
Here $j$ is selected sufficiently large after fixing $L$, $k$ depends on $j$, and $a$ depends on $j,k$. The constant $B_L$ is independent of them.

It is easy to see that $\pi^{-1}(\Sigma_f)$ is a proper analytic subset of $\widetilde Y$ containing the exceptional locus of $\pi$, and hence has zero $\nu$-volume. The same holds for $\pi^{-1}(Z)$, since the blowup of a measure-zero set is still measure-zero.
We may choose a compact set
$\widetilde K_s\subset\pi^{-1}(U_f)$, independent of $j,k,a$, such that
$A\Vol(\widetilde Y\setminus\widetilde K_s,\nu)<1/2$.
By \eqref{eq:scalar-weight-density},
\[
 \sup_{\widetilde Y}w
 \leq A\int_{\widetilde K_s}q\,dV_{\omega_{k,a}}
       +A\Vol(\widetilde Y\setminus\widetilde K_s,\nu) \cdot
          \sup_{\widetilde Y}w
 \leq A\int_{\widetilde K_s}q\,dV_{\omega_{k,a}}
       +\tfrac12\sup_{\widetilde Y}w.
\]
Using $q=e^{2\psi} w$ and \eqref{eq:local-scalar-weight} together with $dV_{\omega_{k,a}}\leq c_s\,dV_{\pi^*\omega}$ on
$\widetilde K_s$, we obtain
\begin{equation}\label{eq:compact-fixed-core}
 \sup_L q
 \leq 2Ac_se^{2B_L}
       \int_{\widetilde K_s} q dV_{\pi^*\omega}, \ \  \text{i.e.} \ \  \sup_L |u_{s,j,k}|_0^2
 \leq 2Ac_se^{2B_L}
       \int_{\widetilde K_s} |u_{s,j,k}|_0^2dV_{\pi^*\omega}.
\end{equation}
Notice that both sides are independent of $a$. 

Keep $j$ fixed and let $k\to\infty$. Since $g_{\epsilon_k}(t) \to g(t)$ for all $t$, the maps $u_{s,j,k}$ and their Hom-norms converge locally smoothly over $U_f$ to $u_{s,j}$ and $|\cdot|_{0,j}$, via the bundle identification induced by the biholomorphism $\pi:\pi^{-1}(U_f)\to U_f$.
Setting $K_s=\pi(\widetilde K_s)$ and passing to the limit in
\eqref{eq:compact-fixed-core}  yields
\begin{equation}\label{eq:physical-fixed-core}
\sup_K |u_{s,j}|_{0,j}^2
\leq C_K \int_{K_s} |u_{s,j}|_{0,j}^2 \ddV,
\quad \forall\; \text{compact set } K\subset U_f
\end{equation}

where $C_K$ is independent of $j$. Clearly the integral over $K_s$ is finite and positive, so normalizing $u_{s,j}$ gives $f_{s,j}$ with
$\int_{K_s}|f_{s,j}|_{0,j}^2\ddV=1$, and \eqref{eq:physical-fixed-core} becomes
\begin{equation}\label{eq:normalized-density}
 \sup_K|f_{s,j}|_{0,j}^2\leq C_K, \quad \forall\; \text{compact set } K\subset U_f.
\end{equation}

Finally, since $f_{s,j}$ satisfies the elliptic equation $\bar\partial^{\operatorname{Hom}(Q_s,T_{s-1,j})} f_{s,j}=0$ with uniformly bounded coefficients (by the inductive hypothesis $T_{s-1,j} \to T_{s-1}$), standard Schauder estimates together with \eqref{eq:normalized-density} give uniform $C^\infty_{\mathrm{loc}}$ bounds for the family $\{f_{s,j}\}$. Thus the Arzelà--Ascoli theorem with the diagonal method provides a subsequence of $f_{s,j}$ that converges locally smoothly to a Higgs morphism $f_{s,\infty}$ on $U_f$, with
\[
 \int_{K_s}|f_{s,\infty}|_{0,\infty}^2\ddV
 =\lim_{j\to\infty}\int_{K_s}|f_{s,j}|_{0,j}^2\ddV=1.
\]
Thus $f_{s,\infty}$ is nonzero. By the same extension argument as for the Higgs field in Proposition~\ref{prop:extension}, $f_{s,\infty}$ extends uniquely to a Higgs morphism on $X$. \end{proof}

\begin{lemma}\label{lem:slope-splitting}
Let $S$ be a stable reflexive Higgs sheaf of slope $\mu$, and let
$\mathcal R=\bigoplus_{\nu\leq\mu}\mathcal R_{\nu}$ be a direct sum of reflexive
polystable Higgs sheaves of slopes $\nu$. If there is a nonzero Higgs morphism
$f:S\to\mathcal R$, then $f$ has zero projection to the blocks with $\nu<\mu$, and
$\mathcal R_{\mu}/f(S)$ is isomorphic to a direct sum of some stable factors of $\mathcal R_{\mu}$.
\end{lemma}

\begin{proof}
We use a classical but important fact about the simplicity of stable sheaves:
\newtheorem*{fact*}{\textbf{Fact}}
\begin{fact*}
If there is a nonzero Higgs morphism $f:S\to T$ between stable reflexive Higgs sheaves $S$ and $T$, then $\mu(S) \leq \mu(T)$ and $f$ is an isomorphism if $\mu(S) = \mu(T)$.
\end{fact*}
The proof of the fact is standard, even in the presence of Higgs fields. We refer to \cite[Theorem~3.1.24(b)]{Kobayashi87} and \cite[Corollary~5.7.12]{Kobayashi87} for details.

The first assertion follows immediately from this fact. Now we decompose
\[
 \mathcal R_\mu=\bigoplus_{\alpha=1}^N\mathcal R_\alpha
\]
into stable reflexive summands of slope $\mu$. Choose a summand projection $p_\alpha: \mathcal R \to \mathcal R_\alpha$ for which $p_\alpha f \neq 0$. The same fact shows that $p_\alpha f$ is an isomorphism. Clearly the morphism $(p_\alpha f)^{-1}p_\alpha : \mathcal R \to S$ is a left inverse of $f$. Thus $\mathcal R=f(S)\oplus\ker p_\alpha$ and
\[
{\mathcal R_\mu}/{f(S)}
\simeq\ker (p_\alpha|_{\mathcal R_\mu})
=
\left(\bigoplus_{\beta\neq\alpha}\mathcal R_\beta\right)
\]
is indeed isomorphic to a direct sum of some stable factors of $\mathcal R_{\mu}$.
\end{proof}

\begin{lemma}
\label{lem:quotient-propagation}
Under the inductive hypotheses at stage $s$, suppose that the morphism $f_{s,\infty}$ given in Proposition~\ref{prop:normalized-maps} is fiberwise injective on $U_f$. Then $P_{s,j}\longrightarrow P_{s,\infty}$ and $T_{s,j}\longrightarrow T_{s}$ in $C^{\infty}_\mathrm{loc}$.
\end{lemma}

\begin{proof}
Recall the local $h_0$-unitary gauges $\{a_{j,\alpha}\}$ over $U_\alpha$ from Proposition~\ref{prop:slice-compactness} such that
\[
A_{j,\alpha}
:=a_{j,\alpha}(A(t_j))
\longrightarrow A_{\infty,\alpha}, \quad
\theta_{j,\alpha}
:=a_{j,\alpha}\theta(t_j)a_{j,\alpha}^{-1}
\longrightarrow\theta_{\infty,\alpha}.
\]

Define $P_{s,j,\alpha}:=a_{j,\alpha}P_{s,j}a_{j,\alpha}^{-1}$ over $U_\alpha \cap U_f$, which can be regarded as a matrix-valued function acting on $(U_\alpha \cap U_f) \times\mathbb{C}^R$ under the local trivialization. Consider the comparison endomorphism
\[
 \begin{aligned}
 C_{s,j,\alpha}
 &=P_{s-1,j,\alpha}P_{s-1,\infty,\alpha}
   +(\Id-P_{s-1,j,\alpha})(\Id-P_{s-1,\infty,\alpha})  \\ \Rightarrow 
 C_{s,j,\alpha}^*C_{s,j,\alpha}
 &=\Id-(P_{s-1,j,\alpha}-P_{s-1,\infty,\alpha})^2.
 \end{aligned}
\]
By the inductive hypothesis $P_{s-1,j} \to P_{s-1,\infty}$, $C_{s,j,\alpha}^*C_{s,j,\alpha}$ is positive definite for sufficiently large $j$, so we may set
\[
 W_{s,j,\alpha}=C_{s,j,\alpha}(C_{s,j,\alpha}^*C_{s,j,\alpha})^{-1/2},
\]
which satisfies
\[
W_{s,j,\alpha}^*W_{s,j,\alpha}=\Id, \quad W_{s,j,\alpha}P_{s-1,\infty,\alpha}=P_{s-1,j,\alpha}W_{s,j,\alpha}, \quad
W_{s,j,\alpha} \longrightarrow \Id \ 
\text{in }C^\infty_{\mathrm{loc}}.
\]
In particular, $W_{s,j,\alpha}:
\operatorname{im}(\Id-P_{s-1,\infty,\alpha})
\xrightarrow{\simeq}
\operatorname{im}(\Id-P_{s-1,j,\alpha})
$ is an isomorphism.

Under the local trivialization induced by the earlier one, we also have matrix-valued functions \(f_{s,j,\alpha}\) on \(U_\alpha \cap U_f\) converging to \(f_{s,\infty,\alpha}\). Since \(f_{s,\infty,\alpha}\) is fiberwise injective, \(f_{s,\infty,\alpha}^*f_{s,\infty,\alpha}\) is positive definite, and we obtain local smooth convergence of \((f_{s,j,\alpha}^*f_{s,j,\alpha})^{-1}\). Notice that via the above local trivialization,
\[
\mathrm{im}(f_{s,j,\alpha})=g_j(\mathcal E_s)/g_j(\mathcal E_{s-1}) \simeq g_j(\mathcal E_s) \cap g_j(\mathcal E_{s-1})^{\perp} = \mathrm{im}(P_{s,j,\alpha}-P_{s-1,j,\alpha}).
\]
By a simple linear-algebraic fact, \(f_{s,j,\alpha}(f_{s,j,\alpha}^*f_{s,j,\alpha})^{-1}f_{s,j,\alpha}^*\) is exactly the orthogonal projection onto \(\mathrm{im}(f_{s,j,\alpha})\). Hence
\[
 P_{s,j,\alpha}
 =P_{s-1,j,\alpha}
   +f_{s,j,\alpha}(f_{s,j,\alpha}^*f_{s,j,\alpha})^{-1}f_{s,j,\alpha}^*.
\]
Taking the limit, we define
\begin{equation}\label{projection-comparison}
 P_{s,\infty,\alpha}
 =P_{s-1,\infty,\alpha}
  +f_{s,\infty,\alpha}(f_{s,\infty,\alpha}^*f_{s,\infty,\alpha})^{-1}f_{s,\infty,\alpha}^*.
\end{equation}
The convergence of $P_{s,j,\alpha}$ and transition functions $g_{i,\alpha\beta}$ in \eqref{eq:gauged-transition-functions} guarantees that we can glue $P_{s,\infty,\alpha}$ to get a global $P_{s,\infty}:E_\infty|_{U_f} \rightarrow E_\infty|_{U_f}$. 

Next we explain the convergence of $T_{s,j}$. Under the same local trivialization as above, we obtain $T_{s,j,\alpha}$, which is isomorphic to $\operatorname{im}(\Id-P_{s,j,\alpha})$. Moreover, using the convergence $P_{s,j,\alpha} \to P_{s,\infty,\alpha}$, we may define the isomorphism 
\[
W_{s+1,j,\alpha}:
\operatorname{im}(\Id-P_{s,\infty,\alpha})
\xrightarrow{\simeq}
\operatorname{im}(\Id-P_{s,j,\alpha})
\]
in the same way as $W_{s,j,\alpha}$. To study the convergence, we identify each $T_{s,j}$ with the single bundle $\operatorname{im}(\Id-P_{s,\infty,\alpha})$. It remains to prove the convergence of the induced holomorphic structures and Higgs fields, which have the following form:
\[
\left\{
\begin{aligned}
 \bar \partial_{j,\alpha} &= (W_{s+1,j,\alpha})^{-1} \cdot (\Id-P_{s,j,\alpha})\bar\partial_{A_{j,\alpha}} \cdot W_{s+1,j,\alpha} \\
\phi_{j,\alpha} &= (W_{s+1,j,\alpha})^{-1} \cdot (\Id-P_{s,j,\alpha})\theta_{j,\alpha} \cdot W_{s+1,j,\alpha}.
\end{aligned}
\right.
\]
Clearly, the pairs $(\bar\partial_{j,\alpha},\phi_{j,\alpha})$ converge to
\[
\left\{
\begin{aligned}
\bar\partial_{\infty,\alpha} &= (\Id-P_{s,\infty,\alpha})\bar\partial_{A_{\infty,\alpha}}, \\
\phi_{\infty,\alpha} &= (\Id-P_{s,\infty,\alpha})\theta_{\infty,\alpha}.
\end{aligned}
\right.
\]
Again we may glue $(\bar\partial_{\infty,\alpha},\phi_{\infty,\alpha})$ to obtain a limiting Higgs bundle $T_s$ over $U_f$, which can be extended to all of $X$ to form a reflexive Higgs sheaf. The convergence $T_{s,j}\to T_s$ is understood as above.
\end{proof}

Finally, we are able to complete the inductive step from stage $s$ to stage $s+1$:

\begin{proof}[Proof of Theorem~\ref{thm:main-theorem}]
Proposition~\ref{prop:HN-vector} tells us that the Harder--Narasimhan type of $(\cE_{\infty},\theta_{\infty})$ is exactly
\[
\frac{\Vol(X,\omega)}{2 \pi}\boldsymbol{\mu}=\frac{\Vol(X,\omega)}{2 \pi}\left(
\underbrace{a_1,\dots,a_1}_{r_1},
\underbrace{a_2,\dots,a_2}_{r_2},
\dots,
\underbrace{a_m,\dots,a_m}_{r_m}
\right),
\]
where $r_s=\operatorname{rank}Q_s$. At stage $s$, we have $(\cE_{\infty},\theta_{\infty})
 \simeq \left(\bigoplus_{q=1}^{s-1}Q_q\right)\oplus T_{s-1}$. Subtracting the slopes (counting multiplicities) accounted for by the split summands $Q_1,\dots,Q_{s-1}$ yields the HN type of $T_{s-1}$:
\[
\frac{\Vol(X,\omega)}{2 \pi}\left(
\underbrace{a_s,\dots,a_s}_{r_s},
\dots,
\underbrace{a_m,\dots,a_m}_{r_m}
\right).
\]
In particular
\[
\mu_{\max}(T_{s-1})=\frac{\Vol(X,\omega)}{2 \pi}a_s=\mu_\omega(Q_s).
\]
By Proposition~\ref{prop:normalized-maps}, we obtain a nonzero Higgs morphism $f_{s,\infty}:Q_s \rightarrow T_{s-1}$. We decompose
\[
T_{s-1}=\bigoplus \{ T_{s-1,\nu}: \nu \leq \frac{\Vol(X,\omega)}{2 \pi} a_s \},
\]
where all $T_{s-1,\nu}$ are polystable. From the proof of Lemma~\ref{lem:slope-splitting}, we know that $f_{s,\infty}$ maps isomorphically onto a stable direct summand of $T_{s-1,\mu_{\omega}(Q_s)}$. Hence $f_{s,\infty}$ is fiberwise injective over $U_f$. 

Now we consider the bundle $\mathcal E_\infty/\operatorname{im}(P_{s,\infty})$ over $U_f$. First, $\operatorname{im}(P_{s,j})=g_j(\mathcal E_s)$ is a holomorphic Higgs subbundle of $g_j(\mathcal E)$, and this is equivalent to the following local equations:
\[
(\Id-P_{s,j,\alpha})
\bar\partial_{A_{j,\alpha}}P_{s,j,\alpha}=0, \quad
(\Id-P_{s,j,\alpha})
\theta_{j,\alpha}P_{s,j,\alpha}=0.
\]
Lemma~\ref{lem:quotient-propagation} then yields the limit
\[
(\Id-P_{s,\infty,\alpha})
\bar\partial_{A_{\infty,\alpha}}P_{s,\infty,\alpha}=0, \quad
(\Id-P_{s,\infty,\alpha})
\theta_{\infty,\alpha}P_{s,\infty,\alpha}=0.
\]
Thus $\operatorname{im}(P_{s,\infty})$ is also a holomorphic Higgs subbundle of $\mathcal E_\infty$, and the quotient $\mathcal E_\infty/\operatorname{im}(P_{s,\infty})$ inherits a Higgs bundle structure, which is isomorphic to $T_s$; the construction of $T_s$ is given in the proof of Lemma~\ref{lem:quotient-propagation}.

Furthermore, by \eqref{projection-comparison} it is easy to see that over $U_f$:
\[
f_{s,\infty}(Q_s) \simeq \frac{\im(P_{s,\infty})}{\im(P_{s-1,\infty})}.
\]
Combining the above facts, we obtain:
\[
\frac{T_{s-1}}{f_{s,\infty}(Q_s)}
\simeq
\frac{
\mathcal E_\infty|_{U_f}/\im(P_{s-1,\infty})
}{
\im(P_{s,\infty})/\im(P_{s-1,\infty})
}
\simeq
\frac{\mathcal E_\infty|_{U_f}}{\im(P_{s,\infty})}
\simeq T_s .
\]
Clearly the above isomorphism can be extended to all of $X$.
Applying Lemma~\ref{lem:slope-splitting} with $f_{s,\infty}:Q_s\to T_{s-1}$, we have
\[
T_{s-1} \simeq f_{s,\infty}(Q_s) \oplus \frac{T_{s-1}}{f_{s,\infty}(Q_s)} \simeq Q_s\oplus T_s,
\]
which establishes inductive hypothesis $(1)$ of stage $s+1$. Inductive hypotheses $(2)$ and $(3)$ follow immediately from Lemma~\ref{lem:quotient-propagation}. At stage $1$, we have $P_{0,j}=0$ and $T_0=\mathcal E_\infty$, so the inductive hypotheses hold automatically, and the process must terminate at stage $m+1$ since $\operatorname{rank} (T_m)=0$, which implies $T_m=0$.
\end{proof}

\section*{Declaration of Generative AI Use}
\addcontentsline{toc}{section}{Declaration of Generative AI Use}

During the preparation of this manuscript, ChatGPT-6 Pro (chat mode) was used to suggest proof strategies for Proposition~\ref{prop:HN-vector} and Proposition~\ref{prop:normalized-maps}, and ChatGPT Codex was used to assist with polishing the English exposition.

The authors have carefully reviewed and verified all mathematical arguments in the manuscript and take full responsibility for its accuracy, completeness, and originality.

\clearpage
\phantomsection

\end{document}